\documentclass[11pt]{amsart}
\usepackage{graphicx}
\usepackage{amsmath}
\usepackage{amscd}
\usepackage{amssymb}
\usepackage{latexsym}
\usepackage{epstopdf}
\newcommand{\R}{\mathbb{R}}
\newcommand{\Q}{\mathbb{Q}}
\newcommand{\C}{\mathbb{C}}
\newcommand{\Z}{\mathbb{Z}}
\newcommand{\Proj}{\mathbb{P}}

\newcommand{\set}[1]{\{ #1 \}}

\newcommand{\norm}[1]{|| #1 ||}
\newcommand{\cl}[1]{\overline{ #1 }}
\newcommand{\del}{\partial}
\newcommand{\delbar}{\overline{\partial}}

\newcommand{\Bl}[1]{\widetilde{#1}}

\newcommand{\pii}[1] {\frac{ #1 }{2\pi i}}
\newcommand{\Tor}{\mathbf{P}}

\newcommand{\ellipstd}[1]{\frac{#1\theta(\frac{#1}{2\pi i}-z)}{\theta(\frac{#1}{2\pi i})}}
\newcommand{\ellip}[1]{\frac{#1\theta(#1-z)\theta'(0)}{\theta(#1)\theta(-z)}}
\newcommand{\ellipar}[2]{\frac{#1\theta(#1-(-#2+1)z)\theta'(0)}{\theta(#1)\theta(-(-#2+1)z)}}

\newcommand{\thetapar}[2]{\frac{\theta(#1-(-#2+1)z)\theta'(0)}{\theta(#1)\theta(-(-#2+1)z)}}
\newcommand{\thetaparr}[2]{\frac{\theta(#1-(-#2+1)z)}{\theta(#1)\theta(-(-#2+1)z)}}

\newcommand{\ellipinv}[1]{\frac{\theta(#1)\theta(-z)}{#1\theta(#1-z)\theta'(0)}}

\newcommand{\jac}[2]{\frac{\theta(#1-(-#2+1)z)\theta(-z)}{\theta(#1-z)\theta(-(-#2+1)z)}}

\newcommand{\ellnorm}[1]{\frac{\theta(#1-z)}{\theta(#1)}}

\newcommand{\orbnormTh}[2]{\frac{\theta(#1-(-#2+1)z)\frac{\theta'(0)}{2\pi i}}{\theta(#1)\theta(-(-#2+1)z)}}
\newcommand{\orbnormZ}[2]{\frac{\theta(#1-(-#2+1)z)\theta(-z)}{\theta(#1)\theta(-(-#2+1)z)}}
\newcommand{\orbellipar}[3]{\frac{#3\theta(#1-(-#2+1)z)\theta'(0)}{\theta(#1)\theta(-(-#2+1)z)}}

\newtheorem{thm*}{Theorem}
\newtheorem{thm}{Theorem}
\newtheorem{prop}{Proposition}
\newtheorem{lem}{Lemma}
\newtheorem{rmk}{Remark}

\newtheorem{cor*}{Corollary}
\newtheorem{cor}{Corollary}

\title{Equivariant Elliptic Genera}
\author{Robert Waelder}
\address{rwaelder@math.ucla.edu}

\begin{document}
\baselineskip = 14 pt
\parskip = 2 pt

\begin{abstract}
    We introduce the equivariant elliptic genus for open varieties and prove an 
    equivariant version of the change of variable formula for 
    blow-ups along complete intersections. In addition, we prove the 
    equivariant elliptic genus analogue of the McKay correspondence 
    for the ALE spaces.
\end{abstract}

\maketitle
\section{Introduction}

The classical McKay correspondence describes a relationship between 
the representation theory of a finite subgroup $G \subset SU(2)$ and 
the topology of the crepant resolution $\Bl{\C^2/G}$ of $\C^2/G$. One 
consequence of this relationship is that the Euler characteristic of 
$\Bl{\C^2/G}$ is equal to the number of irreducible representations 
of $G$. A simple calculation shows that the number of irreducible 
representations of $G$ corresponds in turn to the orbifold Euler number of the 
pair $(\C^2,G)$. Here, if $X$ has an action by a finite group $G$, we define the orbifold Euler 
number $e_{orb}(X,G) =\frac{1}{|G|}\sum_{gh=hg}e(X^{g,h})$, where 
$X^{g,h}$ denotes the common fixed point locus of a pair of commuting 
elements $g$ and $h$. This definition comes from string theory; in 
particular, physicists conjectured that, for $G$ a finite subgroup of
$SU(3)$, the orbifold Euler number 
of $(\C^3,G)$ coincided with the topological Euler number of a crepant 
resolution of the quotient, when such a resolution existed. In analogy 
with the classical McKay correspondence, we refer to 
formulae of this type as McKay correspondences for the Euler 
characteristic.   

Investigations along these lines bring to mind several questions. 
First, what topological data should $e_{orb}(X,G)$ correspond to 
when the quotient $X/G$ does not possess a crepant resolution? Second, what are 
the analogues of the McKay correspondence for other algebro-geometric 
invariants?

In \cite{B}, Batyrev used techniques from motivic integration to define the 
euler number of a pair $(V,D)$ where $D$ is a divisor on $V$. 
The expression $e_{str}(V,D)$ behaves well with 
respect to birational morphisms in the sense that $e_{str}(\Bl{V},\Bl{D}) = 
e_{str}(V,D)$ if $K_{\Bl{V}}+\Bl{D}=\phi^*(K_V+D)$ for a birational 
morphism $\phi: \Bl{V}\rightarrow V$. This definition therefore 
provides a framework for studying the Euler number of a resolution of 
singularities even when no crepant resolution exists. For a special 
choice of a divisor $\Delta$ on $X/G$, Batyrev proved that 
$e_{orb}(X,G) = e_{str}(X/G,\Delta)$. In fact, he proved a much stronger 
variation of this theorem for the $\chi_y$ genus.

An important generalization of both the topological Euler 
characteristic and the $\chi_y$ genus is the two variable elliptic 
genus. If $X$ is an almost complex manifold, the elliptic genus 
$Ell(X)$ is defined as:
$$\int_X 
\prod_{TX}\frac{x_j\theta(\pii{x_j}-z,\tau)}{\theta(\pii{x_j},\tau)}.$$
The product is taken over the formal chern roots of the holomorphic 
tangent bundle to $X$. $\theta(t,\tau)$ is the Jacobi theta 
function and $z$ is a formal parameter. 

When $X$ possesses an action of a finite group $G$, there exists a 
notion of the orbifold elliptic genus of $X$ 
which extends Batyrev's definition of the orbifold $\chi_y$ genus. 
Recently, Borisov and Libgober have proven the elliptic genus analogue 
of the McKay correspondence \cite{BL}. To do this, they first define the 
elliptic genus and orbifold elliptic genus of a pair $(X,D)$ for $D$ a 
divisor on $X$ and show that these definitions satisfy change of 
variable formulae similar to the objects $e_{str}(X,D)$ in Batyrev's paper. 
Whereas Batyrev's proof relies on the change of variables formula from 
motivic integration, Borisov and Libgober examine the case of a single 
blow-up and appeal to the deep result of Wlodarczyk \cite{W} that 
every birational map of smooth complex varieties may be factored into a 
sequence of blow-ups and blow-downs along smooth centers. This change 
of variable formula allowed Borisov and Libgober to reduce the proof 
to a version of the McKay correspondence for toroidal morphisms. A 
crucial aspect of their proof is a description of the cohomological 
pushforward of a toroidal morphism in terms of combinatorial data 
associated to the map. We refer to this technique as Borisov and 
Libgober's push-forward formula.

When $X$ has an action of a compact torus $T$ which commutes with the 
action of a finite group $G$, one has natural definitions for the 
equivariant orbifold elliptic genus of $(X,G)$ and the 
equivariant elliptic genus of $T$-resolutions of $X/G$. One reason 
for studying these equivariant elliptic genera is that, by 
localization, they make sense even when $X$ is not compact, provided 
that $X$ has compact fixed components. In this paper, we prove an 
equivariant elliptic genus analogue of the classical McKay 
correspondence for ALE spaces. Along the way we prove the 
equivariant version of the change of variable formula for blow-ups 
along complete intersections. In the toric case, this change of 
variable formula turns out to be linked to a rigidity property of the 
elliptic genus of a pair. A prominent feature throughout this paper is 
the equivariant analogue of Borisov and Libgober's push-forward formula. 
In this paper we will describe a relationship between Borisov and Libgober's push-forward 
formula and the functorial 
localization formula for a 
toric morphism. The suggestion that Borisov and Libgober's 
push-forward formula is really functorial localization in disguise 
might explain the ease with which the proof of their formula extends 
to the equivariant case.

The sections of this paper are divided as follows: In section \ref{EquivOrbEllClass} we 
introduce the notion of the equivariant orbifold elliptic class, which is useful for making sense of elliptic genera on open varieties with a torus action. In 
section \ref{PrelimEquiv} we discuss various aspects of equivariant cohomology 
and prove some technical lemmas which will be used implicitly 
throughout the paper. In sections \ref{Toric} and \ref{Push} we prove the 
equivariant analogue of Borisov and Libgober's push-forward formula, 
and discuss its relationship to the functorial localization formula 
applied to a toric morphism. In section \ref{Rigid} we prove a rigidity theorem for the elliptic genus of a toric pair $(X,D)$ and discuss its relationship to the change of variable formula, which we prove in section \ref{OrbChVar}.
In section \ref{ALE} we use this results to prove 
an equivariant elliptic genus analogue of the McKay correspondence 
for ALE spaces. Finally, in section \ref{StringyEuler} we discuss the relationship 
between the equivariant elliptic genus and Batyrev's stringy Euler 
number.

\section{Equivariant Orbifold Elliptic Class}\label{EquivOrbEllClass}    
  We begin by setting up some notation.
  Let $X^n$ be a smooth compact variety and $D=\sum_i\delta_iD_i$ a smooth 
  normal crossing divisor, with coefficients $\delta_i<1$. Let $G$ 
  be a finite group acting holomorphically on $X$. For $g,h \in G$ a 
  commuting pair, let $\set{X^{g,h}_\gamma}$ denote the connected components of their
  common fixed point locus. Fix one such component $X^{g,h}_\gamma$.
  The normal bundle $N_{X^{g,h}_\gamma}$ splits as a sum $\oplus_{\lambda} 
  N_{\lambda}$ over irreducible characters for the subgroup $(g,h)$. 
  For $x \in (g,h)$, let $\lambda(x) \in \Q\cap [0,1)$ be the rational 
  number such that $x$ acts on the fibers of $N_{\lambda}$ as 
  multiplication by $e^{2\pi i\lambda(x)}$. 
  
  Now fix an irreducible component $D_i$ of $D$. If $X^{g,h}_\gamma\subset 
  D_i$ then $x\in (g,h)$ acts on the fibers of
  $\mathcal{O}(D_i)|_{X^{g,h}_\gamma}$ as multiplication by $e^{2\pi 
  i\epsilon_i(x)}$ for some rational number $\epsilon_i(x)\in \Q\cap 
  [0,1)$. If $X^{g,h}_\gamma$ is not contained in $D_i$, we define 
  $\epsilon_i = 0$. Of course the functions $\lambda$ and $\epsilon_i$ 
  depend on the choice of the commuting pair $(g,h)$ and on the
  connected component $X^{g,h}_\gamma$ of $X^{g,h}$. We will 
  omit making explicit reference to this dependence in order to simplify the 
  notation.
  
  Following \cite{BL}, we define the orbifold elliptic genus of the pair 
  $(X,D)$ by the formula: 
  
  $$Ell_{orb}(X,D,G)=\frac{1}{|G|}\sum_{gh=hg,\gamma}\int_{X^{g,h}_\gamma}
  \prod_{TX^{g,h}_\gamma}\ellipstd{x_j}\times$$
  $$\prod_{N_{\lambda}}\ellnorm{\frac{x_\lambda}{2\pi 
  i}+\lambda(g)-\lambda(h)\tau}e^{2\pi i\lambda(h)z}\times$$
  $$\prod_{D_i}\jac{\frac{D_i}{2\pi 
  i}+\epsilon_i(g)-\epsilon_i(h)\tau}{\delta_i}e^{-2\pi 
  i\delta_i\epsilon_i(h)z}
  $$
  Here of course $x_j$ denote the Chern roots of $TX^{g,h}_\gamma$, $x_\lambda$ the Chern roots of $N_\lambda$, and $D_i$ the first Chern classes of the corresponding divisors.
  
  Now assume that $X$ has a $T$-action which commutes with the action 
  of $G$, and that the irreducible components of $D$ are $T$-invariant. 
  Assume that the action lifts to the bundles $\mathcal{O}(D_i)$. We 
  define the equivariant orbifold elliptic genus $Ell^T_{orb}$ as follows: For 
  each fixed component $P\subset X^{g,h}_\gamma$, let $\nu_j \in 
  \mathfrak{t}^*$ denote the infinitesimal weights of the torus action 
  on the fibers of the normal bundle $\nu_{P/X^{g,h}_\gamma}$. Similarly, let 
  $\chi_{\lambda} \in \mathfrak{t}^*$ denote the infinitesimal weights of 
  the torus action on $N_{X^{g,h}_\gamma}|_P$. If $P \subset D_i$, let 
  $\eta_i$ denote the infitesimal weight of the torus action on 
  $\mathcal{O}(D_i)|_P$. Otherwise, let $\eta_i = 0$. All 
  of the above weights depend on the fixed component $P$ and
  on the commuting pair $g,h$. Again, we leave 
  this dependence out of the notation in order to avoid cluttering.
  With this in mind, we define $Ell_{orb}^T(X,D,G)=$:
  $$\frac{1}{|G|}\sum_{gh=hg,\gamma}\sum_{P\subset X^{g,h}_\gamma}
  \int_P\prod_{TP}\ellipstd{p_k}\prod_{\nu_P}\ellnorm{\frac{n_j}{2\pi 
  i}+\nu_j}\times$$
  $$\prod_{N_{\lambda}}\ellnorm{\frac{x_\lambda}{2\pi 
  i}+\chi_{\lambda}+\lambda(g)-\lambda(h)\tau}e^{2\pi i\lambda(h)z}\times$$
  $$\prod_{D_i}\jac{\frac{D_i}{2\pi 
  i}+\eta_i+\epsilon_i(g)-\epsilon_i(h)\tau}{\delta_i}e^{-2\pi 
  i\delta_i\epsilon_i(h)z}$$
  
  Finally, motivated by \cite{BL}, we introduce the notion of the equivariant 
  elliptic class $\mathcal{E}ll_{orb}^T(X,D,G) \in H^*_T(X)$. For 
  convenience, assume that every component $X^{g,h}_\gamma$ of $X^{g,h}$ is a 
  connected component of $D_{i_1}\cap\ldots\cap D_{i_r}$ for some 
  indexing set $I^{g,h}_\gamma=\set{i_k}$. We also assume that $D$ is $G$-normal. Then 
  $TX-\oplus_{k=1}^r\mathcal{O}(D_{i_k})\in K_T(X)$ is a bundle which 
  equals $TX^{g,h}_\gamma$ when restricted to $X^{g,h}_\gamma$. Thus, consider 
  the class:
  $$\bigg(\frac{2\pi i\theta(-z)}{\theta'(0)}\bigg)^{n-r}
  \Phi^T_{X^{g,h}_\gamma}\prod_{TX}\ellip{\frac{x_j(t)}{2\pi 
  i}}\prod_{I^{g,h}_\gamma}\ellipinv{\frac{D_i(t)}{2\pi i}}\times$$
  $$\prod_{I^{g,h}_\gamma}\orbnormZ{\frac{D_i(t)}{2\pi 
  i}+\epsilon_i(g)-\epsilon_i(h)\tau}{\delta_i}e^{2\pi i(-\delta_i+1)\epsilon_i(h)z}\times$$
  $$\prod_{I_X-I^{g,h}_\gamma}\jac{\frac{D_k(t)}{2\pi i}}{\delta_k}$$
  Here the parameter $t$ in $x_j(t)$, $D_k(t)$, etc., refers to the equivariant chern 
  roots. When $X$ is simply connected we can always make sense of this definition (see section \ref{EquivChernDiv}). The term $\Phi^T_{X^{g,h}_\gamma}$ is the equivariant Thom class of 
  $X^{g,h}_\gamma$ in $X$. All of the fractions $\epsilon_i$ 
  implicitly depend on $\gamma$. By localization, the integral of this 
  class over $X$ is the 
  contribution from $X^{g,h}_\gamma$ to the equivariant orbifold elliptic genus. We 
  call this equivariant cohomology class 
  $\mathcal{E}ll_{orb}^T(D,X^{g,h}_\gamma)$. We define the equivariant orbifold 
  elliptic class 
  $$\mathcal{E}ll_{orb}^T(X,D,G)=\frac{1}{|G|}\sum_{gh=hg,\gamma}
  \mathcal{E}ll_{orb}^T(D,X^{g,h}_\gamma).$$
  When $G = 1$, we simply call the above expression the equivariant elliptic class $\mathcal{E}ll(X,D)$. We call the integral of this class the equivariant elliptic genus $Ell_T(X,D)$.
     
Suppose that $f:\Bl{X}\rightarrow X$ is the blow-up of $X$ along a $T\times G$-invariant subvariety. Define $\Bl{D}$ on $\Bl{X}$ so that
$f^*(K_X+D) = K_{\Bl{X}}+\Bl{D}$. Then
\begin{thm}
$$f_*\mathcal{E}ll_{orb}^T(\Bl{D},\Bl{X},G) = \mathcal{E}ll_{orb}^T(D,X,G).$$
\end{thm}
This is the equivariant analogue of the change of variable formulae discovered by Chin-Lung Wang and Borisov and Libgober. We will refer to the above formula as the change of variable formula for the orbifold elliptic class. We refer to the formula obtained by integrating both sides as the change of variable formula for the orbifold elliptic genus. For a proof of the formula in this general case, see \cite{RW}. For the purposes of this paper, we will only need to examine the simpler situation in which the blow-up locus and orbifold fixed data of $X$ are complete intersections. For that case, we will provide a complete proof in section \ref{OrbChVar}.

The value of the change of variable formula is that it allows us to compare orbifold elliptic data between varieties which are birationally equivalent. In section \ref{ALE} we will provide an interesting application of the change of variable formula to the computation of equivariant elliptic indices of ALE spaces. The approach we take is inspired by Borisov and Libgober's proof of the non-equivariant McKay correspondence for the elliptic genus.

\begin{rmk}\rm
Note that in case $X$ has a torus action with compact fixed components, via the localization formula we can always make sense of the quantity $Ell_{orb}^T(X,D,G)$ even when $X$ is \it{open}\rm. We will continue to refer to this quantity as the equivariant orbifold elliptic genus of $X$.
\end{rmk}\rm
\begin{rmk}\rm A word on notation: There are many objects associated 
to a variety $X$ which encode the data of the equivariant elliptic genus of $X$. We use 
the prefix $Ell_T$ to refer to objects in $H^*_T(pt)$, $\mathcal{E}ll_T$ 
to refer to objects in $H^*_T(X)$, and $\mathcal{ELL}_T$ to refer to 
objects in $K_T(X)$. 
\end{rmk}

\section{Preliminaries on Equivariant Cohomology}\label{PrelimEquiv}
In this section we gather the ingredients from equivariant 
cohomolgy which we will be using in this paper. For a thorough 
reference on the subject, see \cite{AB}.
\subsection{Definitions and Localization}
Let $X$ be a smooth $T$-space, where $T$ is a compact torus of
rank $\ell$. Let $ET = (S^{\infty})^\ell$. $ET$ is a contractible 
space on which $T$ acts freely. The diagonal action of $T$ on $X\times 
ET$ therefore gives rise to a smooth (infinite-dimensional) quotient 
$X_T = (X 
\times ET)/T$. It is easy to see that $X_T$ 
is a fiber bundle over $BT = ET/T$ with fiber $X$. Define the 
equivariant cohomology group $H^*_T(X) = H^*(X_T)$. 

The translation of concepts from cohomology to equivariant cohomology 
is more or less routine. For example, a $T$-map $f: X\rightarrow Y$ 
gives rise to a natural map $f_T : X_T\rightarrow Y_T$, and therefore 
induces a pullback $f^*:H^*_T(Y)\rightarrow H^*_T(X)$. Similarly, for 
any $E \in K_T(X)$, $E_T$ defines a finite rank vectorbundle over $X_T$ which 
corresponds to the vectorbundle $E\rightarrow X$ over every fiber of 
$X_T \rightarrow BT$. In this way, we may define the equivariant 
characteristic classes of $E$ to be the characteristic classes of 
$E_T$. 

If $p$ is a single point with trivial $T$-action, the equivariant 
map $\pi : X \rightarrow p$ induces a map $\pi^*:H^*_T(p)\rightarrow 
H^*_T(X)$. Since $H^*_T(p) = H^*(BT) = \C[u_1,\ldots,u_\ell]$, the 
map $\pi^*$ makes $H^*_T(X)$ into a $\C[u_1,\ldots,u_\ell]$-module. 
Define $H^*_T(X)_{loc} = 
H^*_T(X)\otimes_{\C[u_1,\ldots,u_\ell]}\C(u_1,\ldots,u_\ell)$. A 
fundamental result of the subject is the localization theorem:

\begin{thm}
    Let $\set{P}$ denote the set of $T$-fixed components of $X$. Then
    $H^*_T(X)_{loc} \cong \bigoplus_P H^*(P)\otimes 
    \C(u_1,\ldots,u_\ell)$.
\end{thm}

If $P$ is a fixed component of $X$, the normal bundle to $P$ splits as 
a sum over the characters of the $T$-action on the fibers: $N_P = 
\bigoplus_{\lambda}V_{\lambda}$. Let $n^i_{\lambda}$ denote the formal chern roots 
of $V_{\lambda}$. If we identify the equivariant parameters 
$u_1,\ldots,u_\ell$ with linear forms on the Lie algebra of $T$, then
the equivariant Euler class $e(P)$ of $N_P$ is equal to 
$\prod_{\lambda}\prod_{i}(n^i_{\lambda}+\lambda)$. Since 
none of the characters $\lambda$ are equal to zero, we see that $e(P)$ 
is always invertible. In light of this fact, we can describe the above 
isomorphism more explicitly. The map $H^*_T(X)_{loc} \rightarrow 
\bigoplus_P H^*(P)\otimes\C(u_1,\ldots,u_\ell)$ is given by $\omega 
\mapsto \bigoplus_P \frac{i_P^*\omega}{e(P)}$, where $i_P :P
\hookrightarrow X$ is the inclusion map.

If $f : X \rightarrow Y$ is a proper map of $T$-spaces, we have the 
equivariant analogue of the cohomological push-forward $f_*: 
H^*_T(X)\rightarrow H^*_T(Y)$. As in the non-equivariant setting, 
$f_*$ satisfies the projection formula $f_* (f^*(\omega)\wedge\eta) = 
\omega\wedge f_*\eta$. The new feature in equivariant cohomology is 
that we have an explicit expression for the restriction of 
$f_*\omega$ to a fixed component in $Y$. This is given by the 
functorial localization formula \cite{MirrorI} \cite{MirrorII}:

\begin{thm}
    Let $f: X\rightarrow Y$ be a proper map of $T$-spaces. Let $P$ be 
    a fixed component of $Y$ and let $\set{F}$ be the collection 
    of fixed components in $X$ which $f$ maps into $P$. Let $\omega 
    \in H^*_T(X)$. Then:
    
    $$\sum_F f_*\frac{i_F^*\omega}{e(F)} = 
    \frac{i_P^*f_*\omega}{e(P)}.$$
\end{thm}

Suppose $f: X\rightarrow Y$ is a proper map of $n$-dimensional $T$ spaces with 
isolated fixed points. For $F$ a fixed point in $X$, let 
$\lambda(F)_1+\ldots+\lambda(F)_n$ denote the decomposition of $T_FX$ 
into irreducible characters. Clearly $e(F) = \prod_{j=1}^n\lambda(F)_j$.
Moreover, each $\omega \in H^*_T(X)$ is defined by a 
collection of polynomial functions (with relations) $\omega_F \in \C[u_1,\ldots,u_\ell]$ 
attached to the fixed points $F$ in $X$. For $P$ a fixed point in $Y$, 
we have by functorial localization:
$$(f_*\omega)_P = \sum_F 
\omega_F\prod_{j=1}^n\frac{\lambda(P)_j}{\lambda(F)_j}.$$ 
In the next two sections we will discuss the similarity between this 
formula and the push-forward formula of Borisov and Libgober.

Note that the localization techniques discussed here continue to hold in the ring formed by uniformly convergent power series of equivariant classes, which is more precisely the domain of definition for the equivariant elliptic class. For simplicity of exposition, we will not make that distinction here. However, see \cite{RW} for a discussion of this technical point.

Before ending this subsection, we make one final remark on an alternative 
approach to equivariant cohomology. Let $e_1,\ldots,e_\ell$ form a 
basis for the Lie algebra of $T$ which is dual to the linear forms 
$u_1,\ldots,u_\ell$. Every $V \in \mathfrak{t}$ defines a vectorfield 
$V$ on $X$ by the formula $V(p) = 
\frac{d}{dt}|_{t=0}\mathrm{exp}(tV)\cdot p$. Define $\Omega^*_T(X)$ to 
be the ring of differential forms on $X$ which are annihilated by 
$\mathcal{L}_V$ for every $V \in \mathfrak{t}$. If we let 
$d_{\mathfrak{t}} = d+\sum_{\alpha=1}^{\ell}u_\alpha i_{e_\alpha}$, 
then $d_{\mathfrak{t}}$ defines an operator on $\Omega^*_T(X)\otimes \C[u_1,\ldots,u_\ell]$ and 
satisfies $d_{\mathfrak{t}}^2 = 0$. The Cartan model for equivariant 
cohomology is defined to be:
$$H^*_T(X)_{\mathrm{Cartan}} = \frac{\ker d_{\mathfrak{t}}}
{\mathrm{im} d_{\mathfrak{t}}}.$$
It is well known that $H^*_T(X)_{\mathrm{Cartan}}\cong H^*_T(X)$. 
See \cite{AB} for details. Throughout, we will switch freely between the two 
descriptions.

\subsection{Equivariant Chern Class of a Divisor}\label{EquivChernDiv}

Let $X$ be a smooth compact simply connected complex manifold with a $T$-action. 
For simplicity of notation, assume that $T = S^1$. We also omit the 
equivariant parameters in this section, since they clutter the 
notation and play no role in the proofs. 
Let $D \subset X$ be a $T$-invariant irreducible Cartier divisor with 
associated line bundle $\mathcal{O}(D)$. Let $\omega$ be a representative of the 
Thom class of the normal bundle, $N_D$, of $D$. By averaging over $T$, we may 
assume that $\omega$ is $T$-invariant. 

Let $V$ be the vectorfield on $X$ induced by the 
$T$-action. We are presented with two 
natural proceedures for extending $\omega$ to an equivariant 
cohomology class, i.e., a $T$-invariant class in the kernel of $d+i_V$. 
First, since $\omega$ is invariant and closed, $\mathcal{L}_V\omega = 
di_V\omega= 0$, so $i_V\omega$ defines a class in $H^1(X)$. Since $X$ 
is simply connected, $i_V\omega = df$. If we require that $f$ have 
compact support in $N_D$, then the above moment map equation defines 
$f$ uniquely, and $\omega - f$ defines an equivariant extension of 
$\omega$.
Second, since $X$ is simply connected and $T$ is abelian, we may lift the 
action of $T$ to $\mathcal{O}(D)$. The equivariant first chern class 
of $\mathcal{O}(D)$ then defines another equivariant extension of $\omega$. 

In this section, we show that both extensions represent 
the same equivariant class provided we choose an ``appropriate'' lift of 
the action of $T$ to $\mathcal{O}(D)$. By appropriate, we mean that 
the action of $T$ on $\mathcal{O}(D)$ extends the natural action of 
$T$ on $\mathcal{O}(D)|_p$ for any fixed point $p$. Note that for dimensionality reasons, the 
two equivariant extensions can differ by at most a constant. The goal 
in this section is to prove that this constant is zero.

\begin{lem}
    Let $\omega$ and $f$ be defined as above. Let $p \in D$ be a fixed 
    point of the $T$-action and let $a=a(p)$ be the infinitesimal weight of 
    the character $\mathcal{O}(D)|_p$. Then $f(p) = -a$.
\end{lem}

\begin{proof}
    Let $U=\set{(z_1,\ldots, z_n)}$ be a coordinate system 
    centered at $p$, with $D$ defined by $\set{z_n = 0}$. We can choose this 
    coordinate system so that $e^{it}\cdot(z_1,\ldots,z_n) = 
    (e^{im_1t}z_1,\ldots,e^{im_{n-1}t}z_{n-1},e^{iat}z_n)$. Let 
    $d\theta$ be the $T$-invariant angular form in the $z_n$-coordinate 
    plane. Call this plane $U_n$. Let $r$ be the distance function on 
    $U_n$ and $\rho(r)$ a bump function which integrates to $1$ over 
    $U_n$ and is identically equal to $1$ in a neighborhood of the 
    origin. Then the Thom class corresponding to the hyperplane $z_n=0$ is 
    represented by the $T$-invariant form $d(\rho(r)d\theta)$ in this 
    neighborhood. It follows that $\omega|_U = d(\rho(r)d\theta) + 
    d\psi$, where $\psi$ is a form with compact support in the 
    $z_n$-direction. Since all the forms involved are $T$-invariant, 
    we may assume that $\psi$ is $T$-invariant. In this coordinate 
    system, the vectorfield
    $V$ takes the form 
    $a\frac{\partial}{\partial\theta}$ in the $U_n$ plane. Thus $i_V 
    d(\rho d\theta) = -d(i_V\rho d\theta) = -d(a\rho(r))$. We therefore 
    have that $i_V\omega = -d(a\rho) - d(i_V\psi)$. Since $a\rho + 
    i_V\psi$ have compact support in the vertical direction and 
    satisfy $d(a\rho + i_V\psi) = -df$, we must have $f = -a\rho-i_V\psi$. 
    Since $V(p) = 0$, this implies that $f(p) = -a$. 
\end{proof}

We next prove that we can always adjust the action of $T$ on 
$\mathcal{O}(D)$ so that it coincides with the natural action of 
$\mathcal{O}(D)|_p$ for any fixed point $p \in D$. The ensuing 
discussion follows closely the ideas of section $8$ of \cite{AB}. Let $\nabla$ be 
a $T$-invariant connection on $\mathcal{O}(D)$ with corresponding 
connection $1$-form $\theta$. View the vectorfield 
$V$ as an operator acting on $\Gamma(X,\mathcal{O}(D))$. If $s$ is a 
local frame, then $Vs = L(s)s$ for some smooth 
function $L$ which depends on $s$. The statement that $\nabla$ is 
$T$-invariant means that $\nabla V = V\nabla$. Fix a local frame $s$ 
satisfying $ds = 0$ in local coordinates. Then $\nabla Vs = \nabla 
Ls = dLs + L\theta s$ and $V\nabla s = V\theta s = \mathcal{L}_V\theta 
s + \theta Ls$. It follows that 
$$dL(s) = \mathcal{L}_V\theta = i_V 
d\theta + di_V\theta.$$ 

Now $d\theta = -2\pi i \omega + d\eta$ for 
some $T$-invariant $1$-form $\eta$. Thus, $i_V d\theta = -2\pi i 
df-d(i_V\eta)$. Thus, the above equation implies $L(s) = 
-2\pi i f+i_V\theta-i_V\eta+2\pi ic$, where $c$ is a constant. It is easy to check that 
this constant is independent of the section $s$. It follows that the 
infinitesimal action of $T$ on $\mathcal{O}(D)|_p$ is given by $-2\pi 
if(p)+2\pi ic$. Thus, the infinitesimal weight attached to every fixed point $p$
is $a(p)+c$. If $c \neq 0$, we can replace $\mathcal{O}(D)$ with 
$\mathcal{O}(D)\otimes\mathcal{O}_c$, where the action of $T$ on 
$\mathcal{O}_c$ takes the global section $1$ to $e^{-2\pi ic}$. Thus, 
we have proven that we can always lift the action of $T$ so that it 
coincides with the natural action on $\mathcal{O}(D)|_p$ for fixed 
points $p$. Whenever we speak of $\mathcal{O}(D)$ as an equivariant 
bundle, we will assume this choice of a lifted action.

Finally, we prove that the equivariant first chern class of 
$\mathcal{O}(D)$ coincides with $\omega-f$. By localization, it 
suffices to prove that $c_1^T(\mathcal{O}(D)) = \omega-f$ at every 
fixed point. But this follows from the well-known observation that 
$c_1^T(\mathcal{O}(D))|_p = a(p)$.

\section{Toric Varieties and Equivariant Cohomology}\label{Toric}
For a good reference on toric varieties, see \cite{F}.
Let $X$ be a smooth complete toric variety of dimension $n$. We denote the fan of $X$ 
by $\Sigma_X$, the lattice of $X$ by $N_X$, and the big torus by $T_X$. Let $Y$ be a smooth 
complete toric variety which satisfies the following properties:

$(1)$: $N_X \subset N_Y$ is a finite index sublattice.

$(2)$: $\Sigma_X$ is a refinement of $\Sigma_Y$ obtained by adding 
finitely-many one dimensional rays.

There is an obvious map of fans $\nu :\Sigma_X \rightarrow \Sigma_Y$ 
which induces a smooth map $\mu : X \rightarrow Y$. We call a map 
induced by such a morphism of fans a \it{toric morphism}\rm. It is 
easy to verify that $\mu: T_X \rightarrow T_Y$ is a covering map with 
covering group $N_Y/N_X$. Thus, we may regard $Y$ as a $T_X$-space. 
Our goal in this section is to obtain a convenient description of the 
equivariant pushforward $\mu_* : H^*_T(X)\rightarrow H^*_T(Y)$ in 
terms of the combinatorics of $\Sigma_X$ and $\Sigma_Y$. Here $T = 
T_X$. 

We first note that fixed points $F$ of $X$ are in $1-1$ 
correspondence with $n$-dimensional cones $C_F \subset \Sigma_X$. 
Furthermore, the infinitesimal weights of the $T$-action on $N_F$ 
correspond to linear forms in $\mathrm{Hom}(N_X,\Z)$ which are dual 
to the generators of $C_F$ in $N_X$. With this in mind, we have the 
following theorem: Let $\C[\Sigma_X]$ denote the ring of piecewise 
polynomial functions on the fan of $X$. Then:

\begin{thm}
    $H^*_T(X) \cong \C[\Sigma_X]$.
\end{thm}

\begin{proof}
    The map $H^*_T(X)\rightarrow \C[\Sigma_X]$ is defined as follows: 
    $\omega \mapsto \set{\omega|_F}_{F\in X^T}$. The fact that the 
    polynomial functions $\omega|_F$ piece together into a 
    well-defined piecewise polynomial function follows from the fact 
    that $\omega$ is a globally defined cohomology class. To define 
    the reverse arrow, it suffices to describe it for piecewise linear 
    functions. If $f \in \C[\Sigma_X]$ is piecewise linear, then it is 
    well-known in toric geometry that $f$ defines a $T$-Cartier 
    divisor $\mathrm{div}(f)$. Let $f \mapsto \mathrm{div}(f)^{\#}$, 
    where $\mathrm{div}(f)^{\#}$ denotes the equivariant extension whose 
    restriction to a fixed point $F$ is $f|_F$.
\end{proof}

Via the identification $H^*_T(X) \cong \C[\Sigma_X]$, we define 
$\nu_*:\C[\Sigma_X]\rightarrow \C[\Sigma_Y]$ to be the map which makes 
the following diagram commute:
$$\begin{CD}
     \C[\Sigma_X]  @>\nu_*>>   \C[\Sigma_Y]\\
    @|                              @|\\
     H^*_T(X)      @>\mu_*>>   H^*_T(Y)\\ 
\end{CD}$$
Here we understand $\C[\Sigma_Y]$ to be the ring of piecewise 
polynomial functions on $\Sigma_Y$ with respect to the lattice $N_X$.

We now describe $\nu_*$ more explicitly. First notice that for $f \in 
\C[\Sigma_X]$, $\nu_*f$ is given by viewing $f|_F$ as the zero degree 
part of an equivariant cohomology class $\omega \in H^*_T(X)$, 
pushing $\omega$ forward by $\mu_*$, and then forming the piecewise 
polynomial function defined by the zero degree part of $\mu_*\omega$. 
Thus, let $C \subset \Sigma_Y$ be an $n$-dimensional cone. Let 
$\nu^{-1}C$ be the fan $\Sigma_C \subset \Sigma_X$ which is the union 
of $n$-dimensional cones $C_i$. Let $x^{C_i}_1,\ldots,x^{C_i}_n$ be 
the linear forms dual to $C_i$ and $x^C_1,\ldots,x^C_n$ the linear 
forms in $\mathrm{Hom}(N_Y,\Z) \subset \mathrm{Hom}(N_X,\Z)$ dual to $C$. 
By functorial localization:
$$(\nu_*f)_C = \sum_{C_i \subset \Sigma_X}f_{C_i}
\frac{\prod_{j=1}^{n} x^C_j}{\prod_{j=1}^{n} x^{C_i}_j}.$$

Similarly, we define $\nu^* : \C[\Sigma_Y]\rightarrow \C[\Sigma_X]$ to 
be the map which makes the following diagram commute:
$$\begin{CD}
     \C[\Sigma_Y]  @>\nu^*>>   \C[\Sigma_X]\\
    @|                              @|\\
     H^*_T(Y)      @>\mu^*>>   H^*_T(X)\\ 
\end{CD}$$

\begin{prop}
    $\nu^*(f) = f\circ \nu$
\end{prop}

\begin{proof}
    Let $\omega \in H^*_T(Y)$ be the form such that $\omega|_P = 
    f|_P$ for every fixed point $P$. Let $F \in \mu^{-1}(P)$. Then
    $$\begin{CD}
     H^*_T(Y)      @>\mu^*>>     H^*_T(X)\\
    @VVV                              @VVV\\
     H^*_T(P)      @>\mu_F^*>>   H^*_T(F)\\ 
\end{CD}$$
commutes. Hence $(\mu^*\omega)|_F = \mu_F^*(\omega|_P) = 
\mu_F^*(f_P) = f_P$. Thus $\nu^*(f)$ is the piecewise polynomial 
function which is equal to $f_{C_P}$ on every cone $C_F \in 
\nu^{-1}C_P$. This is precisely the piecewise polynomial $f\circ \nu$.
\end{proof}

The map $\nu^* : \C[\Sigma_Y] \rightarrow \C[\Sigma_X]$ makes 
$\C[\Sigma_X]$ into a $\C[\Sigma_Y]$-module. As such, we observe:

\begin{prop}
    $\nu_*$ is a $\C[\Sigma_Y]$-module homomorphism.
\end{prop}

\begin{proof}
    In other words, we wish to prove the projection formula 
    $\nu_*(f\nu^*g) = \nu_*(f)\cdot g$. This follows from 
    identifying $\nu_*$ with $\mu_*$, $\nu^*$ with $\mu^*$ and 
    invoking the projection formula from equivariant cohomology.
\end{proof}

\section{Push-Forward Formula for Toroidal Morphisms}\label{Push}

\subsection{Definitions}
Let $X$ be a compact complex manifold and $D_X = \sum_{I_X} D^X_i$ a divisor on $X$ 
whose irreducible components are smooth normal crossing divisors.
For $I\subset I_X$, let $X_{I,j}$ denote the $j$th connected component 
of $\cap_I D^X_i$. Let $X^{o}_{I,j} = X_{I,j}-\cup_{I^c}D^X_i$. The 
collection of subvarieties $X^{o}_{I,j}$ form a stratification of $X$.
Associated to these data is a polyhedral complex with integral 
structure defined as follows:

Corresponding to $X_{I,j}$, define 
$N_{I,j} = \Z e_{i_1,j}+\ldots+\Z e_{i_k,j}$ 
to be the free group on the elements $e_{i_1,j},\ldots,e_{i_k,j}$. 
Here $i_1,\ldots i_k$ are the elements of $I$. 
Define $C_{I,j}$ to be the cone in the first orthant of this lattice. 
Whenever $I'\subset I$ and $X_{I,j}\subset X_{I',j'}$ we have natural 
inclusion maps $N_{I',j'}\hookrightarrow N_{I,j}$ and 
$C_{I',j'}\hookrightarrow C_{I,j}$. Define $\Sigma_X$ to be the 
polyhedral complex with integral structure obtained by gluing the 
cones $C_{I,j}$ together according to these inclusion maps.

Let 
$\C[\Sigma_X]$ denote the ring of piecewise polynomial functions on 
$\Sigma_X$. Fix $C \subset \Sigma_X$. 
Define $f^C$ to be 
the piecewise polynomial function which is equal to 
$\prod_{j=1}^{\dim C}x^{C}_j$ on every cone containing $C$, and 
equal to zero everywhere else. As in the toric geometry case, there is a natural 
correspondence between piecewise linear functions on $\Sigma_X$ and 
Cartier divisors whose irreducible components are components of 
$D_X$. We denote the piecewise linear function corresponding to $D$ 
by $f^D$.

\subsection{Toroidal Morphisms}
Our primary interest in this section is the study of toroidal morphisms.
This is a map $\mu: (X,D_X,\Sigma_X) \rightarrow (Y,D_Y,\Sigma_Y)$ 
which satisfies the following:

$(1)$: $\mu : X-D_X \rightarrow Y-D_Y$ is an unramified cover.

$(2)$: $\mu$ maps the closure of a stratum in $X$ to the closure of a 
stratum in $Y$.

$(3)$: Let $U_y$ be an analytic neighborhood of $y \in Y$ such that the 
components of $D_Y$ passing through $y$ correspond to coordinate 
hyperplanes. Then for $x \in \mu^{-1}(y)$, there exists an analytic 
neighborhood $U_x$ of $x$ such that the components of $D_X$ passing 
through $x$ correspond to coordinate hyperplanes of $U_x$. Moreover, 
the map $U_x \rightarrow U_y$ is given by monomial functions in the 
coordinates. 

Corresponding to $\mu$, we can define a map $\nu : \Sigma_X \rightarrow \Sigma_Y$
as follows: Let $C_{I,i} \subset \Sigma_X$ and let $e_1,\ldots,e_k 
\in N_{I,i}$ be the generators of $C_{I,i}$ which correspond to the 
divisors $D^X_1,\ldots, D^X_k$. We have that 
$\mu(X_{I,i}) = Y_{J,j}$. Let $v_1,\ldots,v_\ell \in N_{J,j}$ be 
the generators of $C_{J,j}$ which correspond to the divisors 
$D^Y_1,\ldots,D^Y_\ell$. For $1 \leq s \leq k$, $1 \leq t \leq 
\ell$, define $a_{st}$ to be the coefficient of $D^X_s$ of the divisor 
$\mu^*(D^Y_t)$. Then we define $\nu(e_s) = \sum a_{st}v_t$. Note that 
if $(X,\Sigma_X) \rightarrow (Y,\Sigma_Y)$ is a smooth toric morphism 
of toric varieties, then $\nu :\Sigma_X \rightarrow \Sigma_Y$ is the 
natural morphism of polyhedral complexes.

We have the following proposition relating $\nu$ to $\mu$:

\begin{prop}\label{Axioms}
    If $C = C_{J,j} \subset \Sigma_Y$, then $\nu^{-1}C$ is the union 
    of fans $\Sigma_\alpha \subset \Sigma_X$ with the following 
    properties:
    
    $(1)$: $\Sigma_\alpha$ is a refinement of $C$ obtained by adding 
    finitely-many $1$-dim rays.
    
    $(2)$: The lattice $N_\alpha$ of $\Sigma_\alpha$ is a finite index 
    sub-lattice of $N_C$.
    
    $(3)$: The fans $\Sigma_\alpha$ are in $1-1$ correspondence with 
    connected components $U_\alpha$ of $\mu^{-1}(N_{Y_{J,j}^{o}})$. 
    The map $U_\alpha \rightarrow N_{Y_{J,j}^{o}}$ is a fibration 
    given by the smooth toric morphism 
    $\Tor_{\Sigma_\alpha,N_\alpha}\rightarrow \Tor_{C,N_C}$ along the 
    fiber, and a $d_\alpha = d(\Sigma_\alpha)$-cover of $Y_{J,j}^{o}$ along the base.
\end{prop}

For a proof, see \cite{BL}. In the examples studied in this paper, it is 
easy to see that the proposition holds. For the purposes of this 
paper, therefore, one may take proposition \ref{Axioms} as an axiom. 

\subsection{Pushforward formula for Polyhedral 
Complexes}\label{Poly}
Motivated by the description of the push-forward $\nu_*$ 
for toric morphisms, define $\nu_* : \C[\Sigma_X] \rightarrow 
\C[\Sigma_Y]$ as follows. Let $C \subset \Sigma_Y$ be an 
$n$-dimensional cone with dual linear forms $x^C_1,\ldots,x^C_n$. 
Then for $f \in \C[\Sigma_X]$, we define:
$$(\nu_*f)_C = \sum_\alpha d_\alpha \sum_{C_i \in \Sigma_\alpha}
f_{C_i}\cdot \frac{\prod_{j=1}^{n}x^C_j}{\prod_{j=1}^{n}x^{C_i}_j}$$
The second sum is taken over the cones $C_i \subset \Sigma_\alpha$ 
with the same dimension as $C$.

Let $V$ be the toric variety $\coprod_\alpha d_\alpha\cdot 
\Tor_{\Sigma_\alpha,N_\alpha}$ with polyhedral fan $\Sigma_V$. We 
have a natural toric morphism $V \rightarrow \C^n$. We can compactify 
$V$ and $\C^n$ to obtain a smooth toric morphism $\overline{V}\rightarrow 
\Proj^n$. If we view $f$ as a piece-wise polynomial function on the 
fan of $\overline{V}$, then the above formula simply corresponds to 
$(\nu_*f)_C$ where $\nu :\Sigma_{\overline{V}} \rightarrow 
\Sigma_{\Proj^n}$. This identification allows us to apply the tools of 
the previous section toward the study of $\nu_*$. 

We first observe that $(\nu_*f)_C$ is indeed a polynomial function. 
This follows from the above identification of $\nu_*$ with the 
equivariant pushforward of a toric morphism. Furthermore, if we 
define $\nu^* :\C[\Sigma_Y]\rightarrow \C[\Sigma_X]$ by the formula 
$\nu^*(f) = f\circ \nu$ then the projection formula:
$$\nu_*(f\nu^*g)_C = \nu_*(f)_C\cdot g_C$$
follows from the projection formula in equivariant cohomology.

\begin{prop}
    $\nu_*(f)$ is a piece-wise polynomial function.
\end{prop}

\begin{proof}
We first show that $\nu_*(f^C)$ is piece-wise polynomial.

Fix $f = f^C$. Suppose $\nu(C) \subset C_0$ for some $C_0 \subset \Sigma_Y$ of 
dimension $k = \dim C$. Then $\nu_*(f)_{C_0} = d(\Sigma_{C_0})\prod_{j=1}^{k}x^{C_0}_j$. 
Suppose $C_1$ is a cone 
containing $C_0$. We wish to show $(\nu_*f)_{C_1}$ is an extension 
of $(\nu_*f)_{C_0}$. 

Consider the toric morphism 
$\sigma:\Tor_{\Sigma_{C_1},N(\Sigma_{C_1})}\rightarrow \C^{\dim C_1}$
induced by the map $\nu: 
\Sigma_{C_1}\rightarrow C_1$. Let $D_1,\ldots, D_k$ be the divisors 
in $\Tor_{\Sigma_{C_1},N(\Sigma_{C_1})}$ which correspond to the 
generators of $C$. Then the piece-wise polynomial 
function $f \in \C[\Sigma_{C_1}]$ represents the equivariant Thom 
class of $D_1 \cap \dots \cap D_k$. Since $\sigma(D_1\cap\dots\cap D_k)$ 
is the affine subspace of $\C^{\dim C_1}$ corresponding to $C_0$, we 
have that $\sigma_*(f)$ is the degree of $\sigma$ along $D_1\cap\dots\cap 
D_k$ times the polynomial function which represents the equivariant 
Thom class of this subspace. But this implies that: 
$$\nu_*(f)_{C_1} =d(\Sigma_{C_1})\frac{[N(\Sigma_{C_1}):N(C_1)]}
{[N(\Sigma_{C_0}):N(C_0)]}\prod_{j=1}^{k}x^{C_0}_j = 
d(\Sigma_{C_0})\prod_{j=1}^{k}x^{C_0}_j.$$

We need to explain the last equality. If $C_0$ corresponds to the 
strata $Y^{o}_{I,j}$ and $U\rightarrow N_{Y^{o}_{I,j}}$ is the fibration 
in Proposition \ref{Axioms} corresponding to the subdivision $\Sigma_{C_0}$, 
then $d(\Sigma_{C_0})[N(\Sigma_{C_0}):N(C_0)]$ and
$d(\Sigma_{C_1})[N(\Sigma_{C_1}):N(C_1)]$ both give the number of 
points in the pre-image of a generic point in $N_{Y^{o}_{I,j}}$. 

Next suppose that $C$ is mapped to a cone $C_0$ of strictly larger dimension. 
Consider the toric morphism 
$\Tor_{\Sigma_{C_0},N(\Sigma_{C_0})}\rightarrow \Tor_{C_0,N(C_0)}$ 
induced by the map $\nu: \Sigma_{C_0}\rightarrow C_0$. The polynomial 
function $f \in \C[\Sigma_{C_0}]$ represents the Thom class of an 
exceptional toric subvariety. Thus $\nu_*(f) = 0$, and it is easy to 
verify that $\nu_*(f) = 0$ on every cone containing $C_0$. Thus, 
$\nu_*$ maps the elements $f^C$ to piecewise polynomial 
functions. Since these functions generate $\C[\Sigma_X]$ as 
a $\C[\Sigma_Y]$-module, the proposition follows from the projection 
formula.
\end{proof}

In what follows we assume that $\mu: X \rightarrow Y$ is an 
equivariant map of compact $T$-spaces. Furthermore, we assume that 
the irreducible components of $D_X$ and $D_Y$ are invariant under 
the $T$-action. Define a map $\rho_X: \C[\Sigma_X]\rightarrow 
H^*_T(X)$ as follows: Fix a cone $C = C_{I,i}$ which corresponds to 
a connected component of the intersection locus of the divisors 
$D_1,\ldots,D_k$. Define $\rho_X[f^C\cdot (f^{D_1})^{a_1}\dots 
(f^{D_k})^{a_k}] = \Phi_{X_{I,i}}\wedge 
D_1^{a_1}\wedge \ldots \wedge D_k^{a_k}$. Here 
$\Phi_{X_{I,i}}$ denotes the (extension by zero) equivariant Thom class of $X_{I,i} \subset X$ 
and, by abuse of notation, $D_j$ denote the (extensions by zero) equivariant Thom classes 
of the divisors $D_j$.

\begin{lem}
    $\rho_X$ is a ring homomorphism.
\end{lem}

\begin{proof}
    Fix cones $C_1=C_{I_1,i_1}$ and $C_2=C_{I_2,i_2}$. It suffices to 
    prove the theorem for the polynomials $f^{C_1}$ and 
    $f^{C_2}$. Let $I=I_1\cup I_2$. Let $C_{I,i}$ denote the cones 
    which correspond to components of the intersection 
    $X_{I_1,i_1}\cap X_{I_2,i_2}$. Clearly 
    $$f^{C_1}f^{C_2} = \sum_{I,i}f^{C_{I,i}}\prod_{I_1\cap I_2}f^{D_j}.$$
    
    Thus $\rho_X(f^{C_1}f^{C_2}) = 
    \sum_{I,i}\Phi_{X_{I,i}}\prod_{I_1\cap I_2}D_j$. However, by the 
    equivariant version of the excess intersection formula, this is 
    precisely the formula for $\rho_X(f^{C_1})\rho_X(f^{C_2})$.
\end{proof}

\begin{lem}\label{CommutPush}
    $\rho_X\nu^* = \mu^*\rho_Y$.
\end{lem}

\begin{proof}
    It suffices to check this for polynomials $f^{C_{I,k}}$.
    If $D$ is a divisor on $Y$ whose irreducible components are 
    components of $D_Y$, then $\nu^*f^D$ is the piecewise linear 
    function corresponding to $\mu^*D$. It follows that 
    $\rho_X\nu^*f^D = \mu^*\rho_Y f^D$. Since all the maps are ring 
    homomorphisms, this implies that $\rho_X\nu^* \prod_{j\in I} 
    f^{D_j} = \mu^*\rho_Y\prod_{j\in I}f^{D_j}$. Let $\mu^*D_i = 
    \sum_j a_{ij}E_j$ as Cartier divisors. As in the lemma in the 
    Appendix, choose equivariant Thom forms $\Phi_{E_j}$ and 
    $\Phi_{D_i}$ with support in small tubular neighborhoods of their 
    respective divisors so that:
    $$\mu^*\Phi_{D_i} = \sum_j a_{ij}\Phi_{E_j} +d\psi_i$$
    as forms. Here $\psi_i$ are equivariant forms with compact support 
    in $\mu^{-1}N_{D_i}$. Let $\set{I,k}$ index the connected 
    components of $\cap_I D_i$. If we choose $N_{D_i}$ sufficiently small, 
    then
    $$\prod_I \Phi_{D_i} = \sum_{I,k}(\prod_I\Phi_{D_i})_{I,k}$$
    where $(\prod_I\Phi_{D_i})_{I,k}$ is the extension by zero of the 
    form $\prod_I\Phi_{D_i}|_{N_{I,k}}$. 
    
    Now $\prod_I f^{D_i} = \sum f^{C_{I,k}}$ and clearly 
    $(\prod_I\Phi_{D_i})_{I,k}$ is a representative of 
    $\rho_Y(f^{C_{I,k}})$. We have that 
    $$\mu^*(\prod_I\Phi_{D_i})_{I,k} = \big\{\prod_I (\sum_j 
    a_{ij}\Phi_{E_j} + d\psi_i) \big\}_{\mu^{-1}N_{I,k}}$$
    where the subscript $\mu^{-1}N_{I,k}$ means the extension by zero 
    of the form restricted to this open set. Since the $\psi_i$ forms 
    have compact support in $\mu^{-1}N_{D_i}$, this form is 
    cohomologous to 
    $$\big\{ \prod_I \sum_j a_{ij}\Phi_{E_j} 
    \big\}_{\mu^{-1}N_{I,k}}.$$
    But this is in turn a representative of $\rho_X\nu^* f^{C_{I,k}}$.
\end{proof}

\begin{lem}\label{PushForward}
    $\mu_* \rho_X = \rho_Y\nu_*$.
\end{lem}

\begin{proof}
    Since $\rho_X\nu^* = \mu^*\rho_Y$ and the polynomials $f^C$ 
    generate $\C[\Sigma_X]$ as a $\C[\Sigma_Y]$-module, by the 
    projection formula it suffices to check $\mu_* \rho_X f^C = 
    \rho_Y\nu_*f^C$.
    
    Case $1$: $C_{I,i}$ is mapped by $\nu$ to a cone $C_{J,j}$ of the same dimension.
    
    From the proof of Proposition $4$, $\nu_*f^{C_{I,i}} = 
    df^{C_{J,j}}$ where $d$ is the degree of $\mu : X_{I,i}\rightarrow 
    Y_{J,j}$. Thus, $\rho_Y\nu_* f^{C_{I,i}} = d\Phi_{Y_{J,j}} = 
    \mu_*\nu_* f^{C_{I,i}}$.
    
    Case $2$: $C_{I,i}$ is mapped by $\nu$ into a cone of strictly 
    larger dimension.
    
    As shown in Proposition $4$, $\nu_*f^{C_{I,i}} = 0$, so 
    $\rho_Y\nu_*f^{C_{I,i}} = 0 = \mu_*\Phi_{X_{I,i}} = \mu_*\rho_X 
    f^{C_{I,i}}$.
\end{proof}

\begin{rmk}\rm
    It is clear that the above lemmas relating $\mu$ to $\nu$ extend 
    without difficulty to the ring $\C[[\Sigma_X]]$ of piecewise 
    convergent power series. 
\end{rmk}

\section{A Rigidity Theorem for Elliptic Genera on Toric Varieties}\label{Rigid}
For $X$ a toric variety and $D \subset X$ a $T$-Cartier divisor, 
the equivariant elliptic genus of the pair $(X,D)$ may be 
interpreted as the equivariant index of an associated differential operator. 
In this section we prove that the equivariant index of this operator 
is actually zero whenever $(X,D)$ satisfies the Calabi-Yau condition $K_X+D=0$. This 
rigidity result closely resembles results by Hattori on the elliptic genera of multifans 
\cite{Hat}.
As we will see, this rigidity theorem is actually closely related to the change of variable
formula for the elliptic genus.
We first define the operator and prove its rigidity.
   
   Let $X$ be a smooth complete toric variety of dimension $n$. Let 
   $T = (S^{1})^{n}$. We can think of $T$ as sitting inside the 
   big-torus of $X$; as such, it induces a natural action on $X$ with 
   isolated fixed points.   
   Let $D_{1},\ldots, D_{\ell}$ be the $T$-invariant divisors 
   corresponding to the one-dimensional cones on the fan of $X$. 
   Suppose $K_{X}+\sum_{i}\delta_{i}D_{i} = 0$ for integers 
   $\delta_{i} \neq 1$. Call such a pair 
   $(X,\sum_{i}\delta_{i}D_{i})$ 
   a toric Calabi-Yau pair. Define 
   $\mathcal{ELL}(\sum_{i}\delta_{i}D_{i})$ to be the following 
   vectorbundle over $X$:
   $$\otimes_{i} \theta'(0)
   \bigotimes_{n=1}^{\infty}
   \Lambda_{-y^{-\delta_{i}+1}q^{n-1}}\mathcal{O}(-D_{i}) 
   \otimes \Lambda_{-y^{\delta_{i}-1}q^{n}}\mathcal{O}(D_{i})\otimes
   S_{q^{n}}\mathcal{O}(-D_{i}) \otimes S_{q^{n}}\mathcal{O}(D_{i}) 
   $$
   The modular properties of the ordinary index of the above operator 
   were discussed by Borisov and Gunnells in \cite{BG}.
   However, they do not prove the rigidity of the 
   equivariant index.
   
   \begin{thm}
       The equivariant index of $\mathcal{ELL}$ is identically zero.
   \end{thm}
   
   \begin{proof}
       We use a modularity argument similar to the one in \cite{L}.
       It suffices to prove that $\mathcal{ELL}$ is rigid under the 
       action of a generic $1$-parameter subgroup $S^{1}\subset T$. We may 
       further assume that this $S^{1}$ action has isolated fixed 
       points. If $p$ is a fixed point of this action, we must have 
       that $p = D_{i_{1}}\cap \ldots \cap D_{i_{n}}$ for some choice 
       of indices $i_{k}$ depending on $p$.  Let $I_{p} = 
       \set{D_{i_{1}},\ldots,D_{i_{n}}}$ and $I_{p}^{c}$ be the 
       remaining $T$-invariant divisors on $X$. $T_{p}X$ splits as:
       $T_{p}X = \mathcal{O}(D_{i_{1}})\oplus \ldots \oplus 
       \mathcal{O}(D_{i_{n}})|_{p}.$ 
       
       Thus, if the exponents of the $S^{1}$ action on 
       $\mathcal{O}(D_{i})$ are $m_{i}$, then the exponents of the 
       action on $T_{p}X$ are $m_{i_{1}},\ldots,m_{i_{n}}$.
       By the fixed point formula, the equivariant index of 
       $\mathcal{ELL}$ is given, up to a normalization factor which is 
       independent of $t$ and the fixed points $\set{p}$, by:
       $$\sum_{p} \prod_{I_{p}}
       \frac{\theta(m_{i}t-(-\delta_{i}+1)z,\tau)}{\theta(m_{i}t,\tau)}
       \prod_{I_{p}^{c}}\theta(-(-\delta_{j}+1)z,\tau).$$
       
       Call this function $F(t,z,\tau)$. Here $\tau \in \mathbb{H}$ 
       is the lattice parameter defining the Jacobi theta function 
       $\theta(t,\tau)$.
       
       Since for $t \in \R$, $F(t,z,\tau)$ is the index of an elliptic 
       operator, have that $F(t,z,\tau)$ is holomorphic for 
       $(t,z,\tau) \in \R \times \C \times \mathbb{H}$. Let us first 
       examine the modular properties of $F$. 
       Define an action of $SL(2,\Z)$ on $\C\times\C\times \mathbb{H}$ 
       by:
       $$\begin{pmatrix} a & b \cr c & d \cr \end{pmatrix} \cdot
       (t,z,\tau) = 
       (\frac{t}{c\tau+d},\frac{z}{c\tau+d},\frac{a\tau+b}{c\tau+d})$$
       
       If $g \in SL(2,\Z)$ and $F$ is a function on 
       $\C\times\C\times\mathbb{H}$, we define $(g\cdot 
       F)(t,z,\tau) = F(g^{-1}(t,z,\tau))$. Let $F$ be the function 
       given by the fixed point formula above. From the relations:
       
       $$\theta \big (\frac{t}{c\tau+d},\frac{a\tau+b}{c\tau+d}\big ) = 
       \zeta(c\tau+d)^{\frac{1}{2}}e^{\frac{\pi 
       ict^{2}}{c\tau+d}}\cdot\theta(t,\tau)$$ 
       we have that 
       $F(\frac{t}{c\tau+d},\frac{z}{c\tau+d},\frac{a\tau+b}{c\tau+d})$ 
       is equal to:
       
       $$\sum_{p}\prod_{I_{p}}\hbox{exp}\big (\frac{\pi 
       ic(m_{i}t-(-\delta_{i}+1)z)^{2}}{c\tau+d}\big ) \cdot 
       \hbox{exp}\big ( \frac{\pi ic(m_{i}t)^{2}}{c\tau+d} \big 
       )^{-1}\cdot$$
       $$\frac{\theta(m_{i}t-(-\delta_{i}+1)z,\tau)}{\theta(m_{i}t,\tau)}
       \prod_{I_{p}^{c}}\hbox{exp}\big (\frac{\pi 
       ic((-\delta_{j}+1)z)^{2}}{c\tau+d} 
       \big)\cdot$$ 
       $$\zeta(c\tau+d)^{\frac{1}{2}} \theta(-(-\delta_{j}+1)z,\tau)$$
       This expression simplifies to:
       $$\zeta^{\ell-n}(c\tau+d)^{\frac{\ell-n}{2}}\hbox{exp}\big (\frac{\pi 
       ic\sum_{i=1}^{\ell}(-\delta_{i}+1)^{2}z^{2}}{c\tau+d}\big )$$       
       $$\sum_{p} \hbox{exp} \big (\frac{-2\pi 
       ic\sum_{I_{p}}m_{i}(-\delta_{i}+1)zt }{c\tau+d}\big )$$
       $$\prod_{I_{p}} \frac{\theta(m_{i}t-(-\delta_{i}+1)z,\tau)}{\theta(m_{i}t,\tau)}
       \prod_{I_{p}^{c}}\theta(-(-\delta_{j}+1)z,\tau).$$
       
       Since $K_{X}+\sum_{i}\delta_{i}D_{i} = 0$ and $K_{X} = 
       -\sum_{i}D_{i}$, we have $\sum_{i}(-\delta_{i}+1)D_{i} = 0$. 
       Thus, the weights at every fixed point for this trivial line 
       bundle must be the same. But the weight at a fixed point $p$ is 
       given by $\sum_{I_{p}}(-\delta_{i}+1)m_{i}$. Since this sum is 
       independent of $p$, we can pull the terms 
       $\hbox{exp}\big (\frac{-2\pi 
       ic\sum_{I_{p}} (-\delta_{i}+1)m_{i}tz}{c\tau+d}\big )$ 
       outside of the summation over the fixed points. We therefore have 
       that $F(g(t,z,\tau)) = K_{g}(t,z,\tau)F(t,z,\tau)$ for some 
       holomorphic nowhere zero function $K_{g}(t,z,\tau)$. In particular, 
       $F(g(t,z,\tau))$ has no poles for $(t,z,\tau) \in 
       \R\times\C\times \mathbb{H}$. 
       
       We now show that $F$ is in fact 
       holomorphic for $(t,z,\tau) \in \C\times\C\times \mathbb{H}$. 
       Clearly, the only poles for $F$ are of the form
       $(\frac{n}{\ell}(c\tau_{0}+d),z_{0},\tau_{0})$, where we may 
       assume that $(c,d) = 1$. Choose integers $a$ and $b$ so that 
       $ad-bc = 1$. Let $g = \begin{pmatrix} a & b\cr c & d \cr 
       \end{pmatrix}$. From the above, we know that $F$ is 
       holomorphic at the point 
       $(\frac{n}{\ell},\frac{z_{0}}{c\tau_{0}+d},
       \frac{a\tau_{0}+b}{c\tau_{0}+d}) = g\cdot 
       (\frac{n}{\ell}(c\tau_{0}+d),z_{0},\tau_{0})$. This implies 
       that $g^{-1}F$ is holomorphic at 
       $(\frac{n}{\ell}(c\tau_{0}+d),z_{0},\tau_{0})$. But 
       $(g^{-1}F)(t,z,\tau) = K_{g^{-1}}(t,z,\tau)F(t,z,\tau)$. Since 
       $K_{g^{-1}}$ is holomorphic and nowhere vanishing, we must have 
       that $F$ is holomorphic at
       $(\frac{n}{\ell}(c\tau_{0}+d),z_{0},\tau_{0})$. Therefore, $F$ 
       is in fact holomorphic on $\C\times\C\times \mathbb{H}$.
       
       Next we prove that $F(t,z,\tau)$ is constant in the 
       variable $t$. Let $z = \frac{1}{N}$ for $N$ an integer. 
       It is easy to verify that 
       $F(t+1,\frac{1}{N},\tau) = F(t,\frac{1}{N},\tau)$ and 
       $F(t+N\tau,\frac{1}{N},\tau) = F(t,\frac{1}{N},\tau)$. Thus, 
       $F(t,\frac{1}{N},\tau)$ is a holomorphic function on a torus, 
       and therefore constant. Hence, for every $N$,
       $$\frac{\partial}{\partial t} F(t,\frac{1}{N},\tau) = 0.$$
       Hence, we must have $\frac{\partial}{\partial t} F(t,z,\tau) = 
       0$. In other words, $F$ is constant in $t$. 

       Finally, since the coefficients $\delta_i \neq 1$, for a generic $S^1$ action
       the summation $\sum_{I_p}(-\delta_i+1)m_i \neq 0$. As in the proof Hattori's vanishing theorems
       for the elliptic genus of multifans \cite{Hat}, we get that $F(t+\tau,z,\tau) =
       e^{2\pi i \sum (-\delta_i+1)m_i z}F(t,z,\tau) = F(t,z,\tau)$, which implies that $F \equiv 0$.
   \end{proof}
   If $\set{D_i}_{i=1}^{\ell}$ are the $T$-Cartier divisors on a toric 
   variety $X$, then $TX$ is stably equivalent to 
   $\bigoplus_{i=1}^{\ell}\mathcal{O}(D_i)$. By the 
   Atiyah-Bott-Lefschetz fixed point formula, this implies that the 
   equivariant index of $\mathcal{ELL}(\sum\delta_iD_i)$ corresponds 
   to the equivariant elliptic genus $\int_X Ell_T(\sum\delta_i D_i)$, 
   up to a normalization factor. 
   
   With the above in mind, we turn our attention to the blow-up of 
   $\C^n$ at the origin. Let $T = (S^1)^n$ act on $\C^n$ as: 
   $(t_1x_1,\ldots,t_nx_n)$. This induces a natural action on $\Bl{\C^n}$. 
   The fixed points of $\Bl{\C^n}$ are the points $p_i = 
   [0:\ldots:1:\ldots:0]$ in the exceptional divisor which have $1$ 
   in the $i$th homogeneous coordinate and zero everywhere else. Set 
   $t_i = e^{2\pi iu_i}$. Then the infinitesimal weights at $p_i$ are 
   $u_1-u_i,\ldots,u_n-u_i,u_i$. 
   
   For $i=1,\ldots,n$, let $\alpha_i < 1$ be the coefficients of the
   coordinate hyperplanes $D_i$.  
   Let $\alpha_0 = \sum_{i=1}^n\alpha_i+(1-n)$. For this simple blow-up, the change
   of variable formula for the equivariant elliptic genus of $(\C^n,\sum\alpha_i D_i)$  takes the following form:
   \begin{lem}\label{BlowUpPt}
       \begin{align*}
	   {}&\sum_{i=1}^{n}\prod_{j\neq i}^n 
	   \thetaparr{u_j-u_i}{\alpha_j}\cdot\thetaparr{u_i}{\alpha_0}\\ 
	   {}&=\prod_{j=1}^n\thetaparr{u_j}{\alpha_j}
       \end{align*}
   \end{lem}
 
   More generally, let $\mu: X\rightarrow Y$ be a composition of toric blow-ups of a smooth complete toric variety $Y$ with associated simplicial map $\nu:\Sigma_X\rightarrow \Sigma_Y$. Since the map has degree $1$, $X$ 
   and $Y$ share the same lattice $N$. For $i=1,\ldots,k$, let $a_i$ denote the 
   $1$-dimensional rays of $\Sigma_X$ and for $j=1,\ldots,\ell$, 
   let $b_j$ denote the $1$-dimensional 
   rays of $\Sigma_Y$. Any sequence $\alpha = 
   \set{\alpha_1,\ldots,\alpha_k}$ of rationals $\alpha_i < 0$ defines a piece-wise linear 
   function $f_\alpha \in \C[\Sigma_X]$ given by $f_\alpha(a_i) = 
   \alpha_i$. This linear function in turn gives rise to the 
   $T$-Cartier divisor $\alpha_1D_{a_1}+\ldots+\alpha_kD_{a_k}$, 
   where $D_{a_i}$ are the divisors associated to the rays $a_i$. 
   
   Clearly $\mu^*(\alpha_1D_{a_1}+\ldots+\alpha_kD_{a_k})$ is the $T$-Cartier 
   divisor on $Y$ corresponding to the linear function 
   $\nu^*f_\alpha$. For each ray $b_j$, let $\beta_j = 
   f_\alpha(b_j)$. Then the sequence $\beta = 
   \set{\beta_1,\ldots,\beta_\ell}$ defines the piece-wise linear 
   function $\nu^*f_\alpha$ and corresponds to the divisor 
   $\mu^*(\alpha_1D_{a_1}+\ldots+\alpha_kD_{a_k})$. We call the 
   sequence $\beta = \mu^*\alpha$ the pull-back of $\alpha$ by $\mu$. 
   
   For each cone $n$-dimensional cone $C_i \subset \Sigma_X$, let $x_{ij}$ denote the 
   linear forms in $\mathrm{Hom}(N,\Z)$ dual to the generators of $C_i$. 
   Let $\beta_{ij}\in \beta$ denote the coefficients corresponding to the 
   generators of $C_i$. For an $n$-dimensional cone $C'_i\subset \Sigma_Y$, define 
   $y_{ij}$ and $\alpha_{ij}$ similarly. Then we have:
   
   \begin{thm}\label{ToricRefine}
       $$\sum_{C_i\subset\Sigma_X}\prod_{j=1}^n
       \frac{\theta(x_{ij}+\beta_{ij}z)}{\theta(x_{ij})\theta(\beta_{ij}z)}=
       \sum_{C'_i\subset\Sigma_Y}\prod_{j=1}^{n}\frac{\theta(y_{ij}+\alpha_{ij}z)}
       {\theta(y_{ij})\theta(\alpha_{ij}z)}.$$
  \end{thm}
  
  \begin{proof}
      By the naturality property of the integer sequences $\alpha$ 
      and $\beta$, if the formula holds for a single blow-up, then it 
      will hold for a composition of blow-ups. Therefore, we can restrict our
attention to the case of a single subdivision.
    Now the divisor $\alpha_1D_1+\ldots+\alpha_kD_k = 
      K_Y+\sum_{i=1}^k(\alpha_i+1)D_{a_i}$. Thus, the right-hand side 
      of the equation in the theorem is just (up to a normalization 
      factor) the equivariant elliptic genus of the pair 
      $(Y,\sum_i(\alpha_i+1)D_{a_i})$. Since 
      $\mu^*(\alpha_1D_1+\ldots+\alpha_kD_k) = 
      \beta_1D_{b_1}+\ldots+\beta_{\ell}D_{b_{\ell}} = 
      K_X+\sum_j(\beta_j+1)D_{b_j}$, the left hand side of the 
      equation is the elliptic genus of $(X,D)$ where $K_X+D = 
      \mu^*(K_Y+\sum_i(\alpha_i+1)D_{a_i})$.

 It clearly suffices to prove that
for $C_i' \subset \Sigma_Y$ an $n$-dimensional cone, the contributions to the RHS coming
from $C_i'$ correspond to the contributions to the LHS coming from the cones $C_i \subset \Sigma_X$
mapping into $C_i'$. Given the above identifications, this amounts to proving the change of variable
formula for the blow-up of $\C^n$ along a $T$-invariant subspace, with the standard torus action.
Since every such blow-up may
be viewed as a product of the identity map along $\C^k$ times the blow-up at the origin of $\C^{n-k}$,
it suffices to prove lemma \ref{BlowUpPt}

Compactify $\C^n$ be viewing it as a subset of $\Proj^n$, and 
       extend the torus action in the obvious manner. We may 
       similarly view $\Bl{\C^n}$ as an open subset of the blow-up 
       of $\Proj^n$ at the origin. Both compactifications 
       are toric varieties and the induced actions are consistent with the 
       action of the big torus. Let $H \subset \Proj^n$ denote the 
       hyperplane at infinity--that is, the hyperplane disjoint from 
       the blow-up point $p_0=[0:\ldots:0:1]$. Since all the divisors $D_i$ 
       corresponding to coordinate hyperplanes passing through $p_0$ 
       are linearly equivalent to $H$ and $K_{\Proj^n} = -(n+1)H$, 
       the line bundle $L=K_{\Proj^n}+\sum_{i=1}^n\alpha_i 
       D_i+((n+1)-\sum_i\alpha_i)H$ is trivial. Thus, the equivariant 
       elliptic genus 
       $\int_{\Proj^n}\mathcal{E}ll_T(\Proj^n,L-K_{\Proj^n})$ is zero. 
       Similarly, $f^*L = 
       K_{\Bl{\Proj^n}}+\sum_{i=1}^n\alpha_i\Bl{D}_i
       +\alpha_0E+((n+1)-\sum_i\alpha_i)f^*H = 0$, which implies 
       that 
       $\int_{\Bl{\Proj^n}}\mathcal{E}ll_T(\Bl{\Proj^n},f^*L-K_{\Bl{\Proj^n}})= 0$. 
       Thus:
       $$\int_{\Bl{\Proj^n}}\mathcal{E}ll_T(\Bl{\Proj^n},f^*L-K_{\Bl{\Proj^n}})=
       \int_{\Proj^n}\mathcal{E}ll_T(\Proj^n,L-K_{\Proj^n}).$$
       It is easy to see that the contribution to the left-hand 
       integral coming from the fixed points mapping to $p_0$ is the LHS of the equation in
       lemma \ref{BlowUpPt}. 
       Similarly, the contribution to the 
       right-hand integral which comes from the blow-up point $p_0$ is 
       equal to the RHS of lemma \ref{BlowUpPt}. Since $\Proj^n$ 
       and $\Bl{\Proj^n}$ are isomorphic away from these points, the 
       contributions to the two integrals coming from the other fixed 
       points are the same, and cancel from both sides of the equation. 
       This proves lemma \ref{BlowUpPt} and completes the proof. 
  \end{proof}

We will see shortly that the addition formula in theorem \ref{ToricRefine} lies at the heart of the change of variable formula.
 
\section{Equivariant Change of Variables Formula}
\subsection{Preliminaries}
Let $(X,D=\sum_{I_X}\alpha_iD_i,G)$ be a $G$-normal pair with $\alpha_i <1$ and with a $T$-action commuting 
      with $G$ and acting invariantly on $D$. Assume that every component
      $X^{g,h}_\gamma$ of $X^{g,h}$ is a complete intersection of 
      components of $D$. Let $f:\Bl{X}\rightarrow X$ be the blow-up of 
      $X$ along a smooth $G$-invariant subvariety whose components are 
      complete intersections of components of $D$. Define 
      $\Bl{D}=\sum_{I_{\Bl{X}}}\delta_j\Bl{D}_j$ so 
      that $K_{\Bl{X}}+\Bl{D}=f^*(K_X+D)$. Note that since 
$\alpha_i < 1$, the coefficient in front of $E$ is less than $1$.  
Our goal in this section is to prove the following equivariant change of variable formula
for the orbifold elliptic genus:
\begin{thm}\label{OrbChVar}
      With the above notation, fix a component $X^{g,h}_\gamma$ and let $\Bl{X}^{g,h}_\mu$ be the 
      components of $\Bl{X}^{g,h}$ which map to $X^{g,h}_\gamma$. Then:
      
      $$f_*\sum_\mu\mathcal{E}ll_{orb}^T(\Bl{D},\Bl{X}^{g,h}_\mu)
      = \mathcal{E}ll_{orb}^T(D,X^{g,h}_\gamma).$$
\end{thm}

Let $I_X$ index the irreducible components of $D$. Let 
      $I_{\Bl{X}}$ index the proper transforms of these components, 
      plus the exceptional divisors. Let $\Sigma_X$ be the polyhedral 
      complex associated to $\set{D_i}_{I_X}$ and let 
      $\Sigma_{\Bl{X}}$ be the polyhedral complex associated to 
      $\set{\Bl{D}_j}_{I_{\Bl{X}}}$. Note that if $X_{I_k,i_k}$ are the 
      components of the blow-up locus, then $\Sigma_{\Bl{X}}$ is 
      obtained from $\Sigma_{X}$ by adding the ray through the point 
      $(1,\ldots,1)$ in each of the cones $C_{I_k,i_k}$. 
The map 
$(\Bl{X},\Sigma_{\Bl{X}},\Bl{D})\rightarrow (X,\Sigma_X,D)$ clearly 
satisfies the axioms of a toroidal morphism.
Before proceeding with the proof, we need to establish some cohomological
properties of this toroidal morphism:

For any variety $X$ with normal crossing divisors $\set{D_i}_{I_X}$, let 
$\Omega(\log D)$ be the locally-free sheaf defined as follows: Let $U = 
\set{(x_1,\ldots,x_k, x_{k+1},\ldots,x_n)}$ be a local coordinate 
system centered at $p$ whose coordinate hyperplanes 
$x_{k+1}=0,\ldots,x_n=0$ correspond to the divisors $D_{k+1},\ldots,D_n \in 
\set{D_i}$ passing through $p$. Then $\Omega(\log D)(U)$ is the 
$\mathcal{O}_X$-module generated 
by the forms $dx_1,\ldots,dx_k, 
\frac{dx_{k+1}}{x_{k+1}},\ldots,\frac{dx_n}{x_n}$. We have the 
following exact sequence of sheaves:
$$0\to \Omega^1 \rightarrow \Omega(\log D) 
\rightarrow \bigoplus \mathcal{O}_{D_i} \to 0.$$
The first map is the obvious inclusion. The second arrow is the residue 
map which takes a section $\omega =\sum_i f_i \frac{dx_i}{x_i}$ to 
$\bigoplus_i f_i|_{D_i}$. It is clear that this map is zero precisely 
when $\omega$ defines a local holomorphic section of $T^*X$. From the 
exact sequence of sheaves
$$0\to \mathcal{O}(-D_i) \rightarrow \mathcal{O}_X \rightarrow 
\mathcal{O}_{D_i}\to 0$$
we get that $\Omega(\log D)-T^*X = -\sum \mathcal{O}(-D_i)$ as 
stable vectorbundles. Applying the dual of this formula to the varieties $(X,D)$ 
and $(\Bl{X},\Bl{D})$ defined above and observing that 
$f^*\Omega(\log D) = \Omega(\log\Bl{D})$, we arrive at the 
K-theoretic relation:
$$T\Bl{X}-f^*TX = 
\sum_{I_{\Bl{X}}}\mathcal{O}(\Bl{D}_j)-\sum_{I_X}f^*\mathcal{O}(D_i).$$
The equality is on the level of stable equivalence. We claim that the 
equality holds in $K_T(\Bl{X})$. To prove this, we verify the equality 
at every fixed component $F \subset \Bl{X}$. Let $F \in f^{-1}(P)$ be 
a fixed component which maps to $P$. Denote by $i_F$ and $i_P$ the
inclusions of $F$ and $P$ in $\Bl{X}$ and $X$. Let $\Bl{D}_1,\ldots,\Bl{D}_\ell$ be the divisors 
on $\Bl{X}$ which contain $F$, and $D_1,\ldots,D_r$ the divisors on 
$X$ containing $P$. Then $i_P^*TX = TP\oplus 
N\oplus\bigoplus_{i=1}^{r}i_P^*\mathcal{O}(D_i)$. Here if $Z$ is the 
connected component of $\cap_{i=1}^{r}D_i$ containing $P$, then $N$ is the 
normal bundle of $P$ in $Z$. Similarly, 
$i_F^*T\Bl{X} = 
TF\oplus\Bl{N}\oplus\bigoplus_{j=1}^{\ell}i_F^*\mathcal{O}(\Bl{D}_j)$, with 
$\Bl{N}$ defined similarly. A computation in coordinates reveals that 
$f^*N = \Bl{N}$. Since $i_F^*f^* = f^*i_P^*$, we have
$$i_F^*T\Bl{X} - i_F^*f^*TX = 
TF+\sum_{j=1}^{\ell}i_F^*\mathcal{O}(\Bl{D}_j)-f^*TP-\sum_{i=1}^{r}i_F^*f^*\mathcal{O}(D_i).$$
However, by the non-equivariant formula for $T\Bl{X}-f^*TX$ derived 
from the log complex, we have that 
$$TF-f^*TP = \sum_{E_j\cap F<F}i_F^*\mathcal{O}(\Bl{D}_j)-
\sum_{D_i\cap P<P}i_F^*f^*\mathcal{O}(D_i)$$
where the sums are taken over the divisors which intersect properly 
with the fixed components. Since these bundles all carry trivial 
$T$-actions, the above formula holds in the equivariant category. 
Finally, observe that if $\Bl{D}_j$ is disjoint from $F$ (resp. $D_i$ is 
disjoint from $P$) then $i_F^*\mathcal{O}(\Bl{D}_j)$ (resp. 
$i_P^*\mathcal{O}(D_i)$) is equivariantly trivial. Hence:
$$i_F^*T\Bl{X} - i_F^*f^*TX = \sum_{I_{\Bl{X}}} 
i_F^*\mathcal{O}(\Bl{D}_j)-\sum_{I_X} i_F^*f^*\mathcal{O}(D_i).$$
From this we deduce the following important formula relating the 
equivariant chern roots of $\Bl{X}$ to $X$:

\begin{lem}\label{ChernRoots}
    Let $f:\Bl{X}\rightarrow X$, $\set{\Bl{D}_j}$, and $\set{D_i}$ be as 
    above. Then:
    
    $$\frac{c_T(T\Bl{X})}{f^*c_T(TX)} = 
    \frac{\prod_{I_{\Bl{X}}}(1+c_1^T(\Bl{D}_j))}{\prod_{I_X}(1+f^*c_1^T(D_i))}.$$
\end{lem}

\subsection{Proof of the Change of Variables Formula}

The method of proof used here is adapted from Borisov and Libgober's 
calculation of the pushforward of the orbifold elliptic genus by a 
toroidal morphism. \cite{BL}

By lemma \ref{ChernRoots}
 $\mathcal{E}ll_{orb}^T(\Bl{D},\Bl{X}^{g,h}_\mu)$ is equal to
    
     $$\Phi^T_{\Bl{X}^{g,h}_\mu}f^*\bigg\{\prod_{TX}\ellip{\frac{x_j(t)}{2\pi 
  i}}\prod_{I_X}\ellipinv{\frac{D_i(t)}{2\pi i}}\bigg\}\cdot$$
  $$\prod_{I^{g,h}_\mu}\thetapar{\frac{\Bl{D}_i(t)}{2\pi 
  i}+\epsilon_i(g)-\epsilon_i(h)\tau}{\delta_i}e^{2\pi i(-\delta_i+1)\epsilon_i(h)z}\cdot$$
  $$\prod_{I_{\Bl{X}}-I^{g,h}_\mu}\ellipar{\frac{\Bl{D}_k(t)}{2\pi i}}{\delta_k}
  \bigg(\frac{2\pi i\theta(-z)}{\theta'(0)}\bigg)^{n}$$
  
  Thus, in order to prove the change of variables formula, we are 
  reduced to proving:
  
  $$f_*\bigg\{\sum_{\Bl{X}^{g,h}_\mu}
  \prod_{I^{g,h}_\mu}
  \orbnormTh{\frac{\Bl{D}_i(t)}{2\pi i}+\epsilon_i(g)-\epsilon_i(h)\tau}
  {\delta_i}e^{2\pi i(-\delta_i+1)\epsilon_i(h)z}\cdot$$
  $$\prod_{I_{\Bl{X}}-I^{g,h}_\mu}\ellipar{\frac{\Bl{D}_i(t)}{2\pi i}}
  {\delta_i}\cdot \Phi^T_{\Bl{X}^{g,h}_\mu}\bigg\} = $$
  
  $$
  \prod_{I^{g,h}_\gamma}
  \orbnormTh{\frac{{D}_i(t)}{2\pi i}+\epsilon_i(g)-\epsilon_i(h)\tau}
  {\alpha_i}e^{2\pi i(-\alpha_i+1)\epsilon_i(h)z}\cdot$$
  $$\prod_{I_{X}-I^{g,h}_\gamma}\ellipar{\frac{D_i(t)}{2\pi i}}
  {\alpha_i}\cdot \Phi^T_{X^{g,h}_\gamma}.$$
  Call the expression in the curly braces $\Omega^{g,h}_\gamma$.
  For each $X^{g,h}_\gamma$, define a piece-wise convergent power 
  series $F^{g,h}_\gamma \in \C[[\Sigma_X]]$ as follows: Let 
  $C^{g,h}_\gamma$ be the cone which corresponds to $X^{g,h}_\gamma$. 
  For $C=C_{I,j}$ a cone containing $C^{g,h}_\gamma$, define:
  $$F^{g,h}_\gamma|_C=
  \prod_{I}\orbellipar{\frac{x^C_i}{2\pi 
  i}+\epsilon_i(g)-\epsilon_i(h)\tau}{\alpha_i}{\frac{x^C_i}{2\pi i}}e^{2\pi 
  i(-\alpha_i+1)\epsilon_i(h)z}$$
  Here for $i\in I$, $x^C_i$ are the linear functions dual to the 
  generators of $C$. If $D_i$ are the divisors which correspond to 
  the generators of $C$, then $\epsilon_i(g)$ and $\epsilon_i(h)$ are the 
  infinitesimal weights of the $g$ and $h$ action on 
  $\mathcal{O}(D_i)|_{X^{g,h}_\gamma}$. Finally, $\alpha_i$ refer to 
  the coefficients of these $D_i$ in the divisor $D$ defined in the 
  statement of the theorem. Finally, for $C$ a cone not containing 
  $C^{g,h}_\gamma$, we define $F^{g,h}_\gamma|_C = 0$. It is easy to see 
  that $F^{g,h}_\gamma$ is a well-defined piece-wise convergent power 
  series and that $\rho_X(F^{g,h}_\gamma)=(\frac{\theta'(0)}{2\pi 
  i\theta(-z)})^n \Omega^{g,h}_\gamma$. We define 
  the piece-wise convergent power series $\Bl{F}^{g,h}_\mu\in 
  \C[[\Sigma_{\Bl{X}}]]$ similarly. By lemma \ref{CommutPush} of section \ref{Poly}, we 
  have reduced the problem to proving:
  $$\nu_*\sum_\mu \Bl{F}^{g,h}_\mu = F^{g,h}_\gamma.$$
  For each cone $C$ containing $C^{g,h}_\gamma$, $\Sigma_C=\nu^{-1}C$ is a 
  subdivision of $C$ obtained by adding no more than one ray through 
  the point $(1,\ldots,1)$ in each subcone. It is 
  clear that the cones $C^{g,h}_\mu$ must be cones inside the fan 
  $\Sigma_C$. Moreover, every cone in $\Sigma_C$ with the same 
  dimension as $C$ will contain exactly one $C^{g,h}_\mu$ as a 
  subcone. (If it contained more than one, that would contradict the 
  fact that the $\Bl{X}^{g,h}_\mu$s are disjoint.) Thus, to prove the 
  formula, we may restrict all our attention to the morphisms 
  $\nu:\Sigma_C\rightarrow C$ for $C \supset C^{g,h}_\gamma$. For 
  $C_j \subset \Sigma_C$ ($\dim C_j=\dim C$), let $x_{ij}=x^{C_j}_i$ and similarly define 
  $\epsilon_{ij}$, $\delta_{ij}$ in the obvious manner. Let $x_i = 
  x^C_i$. By the pushforward formula for $\nu_*$, we are reduced to proving:
  $$\sum_{C_j\subset\Sigma_C}\prod_i\orbnormTh{\frac{x_{ij}}{2\pi 
  i}+\epsilon_{ij}(g)-\epsilon_{ij}(h)\tau}{\delta_{ij}}e^{2\pi 
  i(-\delta_{ij}+1)\epsilon_{ij}(h)z}
  $$
  $$=\prod_i\orbnormTh{\frac{x_{i}}{2\pi 
  i}+\epsilon_{i}(g)-\epsilon_{i}(h)\tau}{\alpha_{i}}e^{2\pi 
  i(-\alpha_{i}+1)\epsilon_{i}(h)z}
  $$
  To prove this, view $\nu: \Sigma_C\rightarrow C$ as a toric 
  morphism. Let $N$ be the lattice corresponding to the toric 
  varieties defined by the fans $\Sigma_C$ and $C$. The elements $g$ 
  and $h$ act on $\Tor_{C,N}$ and $\Tor_{\Sigma_C,N}$ as elements of the big 
  torus. As such, we may view $g$ and $h$ as elements of a 
  sup-lattice $N' \supset N$ of finite index. Under this 
  identification, and using the transformation properties of the Jacobi theta
function, we may assume $\epsilon_i(g) = x_i(g)$ and $\epsilon_{ij}(g) = 
  x_{ij}(g)$. 
  
  Following the notation in the proof of theorem \ref{ToricRefine}, let 
  $\alpha = \set{\alpha_i-1}$ and $\delta = \set{\delta_j-1}$, where the indices 
  range over all one-dimensional rays in $C$ and $\Sigma_C$. 
  Clearly $\delta = \nu^*\alpha$. Let 
  $f_\alpha$ be the linear function on $C$ induced by the 
  multi-index $\alpha$. Then $\sum_i 
  (\alpha_i-1)\epsilon_i(h)=f_\alpha(h)$. Similarly, 
  $\sum_j(\delta_{ij}-1)\epsilon_{ij}(h) = 
  f_\delta(h)=f_{\nu^*\alpha}(h) = \nu^*f_\alpha(h) = f_\alpha(h)$. 
  Thus, all the exponentials on both sides of the above equation are 
  the same. Now the equation follows from Theorem \ref{ToricRefine} after 
  substituting $x_i$ with $x_i+\epsilon_i(g)-\epsilon_i(h)\tau$. This 
  completes the proof.

\subsection{Localization Change of Variable Formula}
Suppose now that $f: \Bl{X}\rightarrow X$ is a $T\times G$ invariant blow-up of open
varieties with compact $T$-fixed components. Assume further that $f:(\Bl{X},\Bl{D})\rightarrow
(X,D)$ admits a smooth equivariant compactification $f':(\Bl{X}',\Bl{D}')\rightarrow (X',D')$
so that the compactified pairs remain $G$-normal, and $f'$ is again an equivariant blow-up. Let
$P$ be a fixed component of $X$ and $\set{F}$ the collection of fixed components mapping to $P$. Then by
the change of variable formula for the orbifold elliptic genus of the compactification, plus functorial localization,
we get the following localized version of theorem \ref{OrbChVar}:
\begin{cor}
$$\sum_{F,\mu}\int_F\frac{i^*_F\mathcal{E}ll_{orb}^T(\Bl{D},\Bl{X}^{g,h}_\mu)}{e(F)}
=\int_P\frac{i^*_P\mathcal{E}ll_{orb}^T(D,{X}^{g,h}_\gamma)}{e(P)}$$
\end{cor}

\section{Equivariant Indices of ALE Spaces}\label{ALE}
  In this section we prove an equivariant elliptic genus analogue of 
  the McKay correspondence for ALE spaces.
  Let $G \subset SU(2)$ be a finite subgroup. Let $T = S^1$ act on 
  $\C^2$ by the diagonal action. $T$ clearly commutes with $G$. Our 
  first goal is to construct an equivariant resolution of 
  singularities for $\C^2/G$.
  
  Let $g,h \in G$ be commuting pairs. Assume one of $g$ or $h$ is 
  nontrivial so that $(\C^2)^{g,h} = (0,0)$. Since $g$ and $h$ 
  commute, they have a simultaneous eigenbasis $(a,b)$ and $(-b,a)$. 
  Let $\ell = \C(a,b)$ and $\ell^{\perp} = \C(-b,a)$. Let 
  $D_{\ell}=\set{-bx+ay=0}$ and $D_{\ell^{\perp}}=\set{ax+by=0}$. 
  Then $D_{\ell}, D_{\ell^{\perp}}$ are $(g,h)\times T$-invariant 
  normal crossing divisors and $D_{\ell}\cap D_{\ell^\perp} = 
  (\C^2)^{g,h}$.
  
  Let $f: \Bl{\C^2}\rightarrow \C^2$ be the blow-up of $\C^2$ at 
  the origin with exceptional divisor $E$. Let $\Bl{D_\ell}$ and 
  $\Bl{D_{\ell^\perp}}$ denote the proper transforms of $D_\ell$ and 
  $D_{\ell^\perp}$.
  The collection of divisors 
  $\set{\Bl{D}_\ell,\Bl{D}_{\ell^\perp}}$ over all commuting pairs 
  $(g,h)$ plus the exceptional divisor $E$ form a system of 
  $T$-invariant normal crossing divisors. Moreover, this system of 
  divisors is $G$-invariant. This just follows from the fact that 
  if $\ell$ is an eigenvector for $g$, then $h\cdot\ell$ is an 
  eigenvector for $hgh^{-1}$. Let $X = \Bl{\C^2}$, and label this 
  system of divisors $\set{D_j}_{I_X}$. The following lemma 
  demonstrates that $(X,\sum D_j,G)$ is, in the words of Batyrev [B], 
  a canonical abelianization of $(\C^2,0,G)$:
  
  \begin{lem}
      $(X,\sum_j D_j,G)$ is $G$-normal with abelian stabilizers.
  \end{lem}
  
  \begin{proof}
      We first prove $G$-normality. 
      Suppose $g\cdot x = x$ and $x \in D_j$, for some $g\neq e$. 
      Then $x = \ell \in E$, since $G$ acts freely on $X-E$. 
      Identifying $\ell$ with a line in $\C^2$, we see that $\ell$ is 
      an eigenspace for $g$, and therefore, that $D_j = 
      \Bl{D}_\ell$. It follows that $g\cdot D_j = D_j$, which proves 
      $G$-normality.
      
      Next, let $g,h$ be two elements which fix $\ell$. Then $\ell$ 
      and $\ell^\perp$ are orthogonal eigenspaces for $g$ and $h$. It 
      follows that $g$ and $h$ are simultaneously diagonalizeable, and 
      therefore commute. This completes the proof.
  \end{proof}
      
  Note that the above proof also implies that $X^g = \cap_J D_j$ for 
  some indexing set $J \subset I_X$.
  
  Consider the map $f: X \rightarrow X/G$. Let $p \in X$ have a 
  nontrivial stabilizer $G_p$. Notice that $p \in X^T$ and $p = 
  D_p \cap E$ for some divisor $D_p \in \set{D_j}_{I_X}$. Let $U_p = 
  \set{(x,y)}$ be a coordinate system centered at $p$ such that 
  $\set{x=0} = E\cap U_p$ and $\set{y=0}=D_p\cap U_p$. Then $f(U_p) 
  \cong U_p/G_p$. Since $(X,\sum_jD_j,G)$ is $G$-normal, we may 
  identify $G_p$ as a finite subgroup of $S^1\times S^1$ acting on 
  $U_p$ in the obvious manner. Thus, let $\Bl{U_p/G_p}$ be a toric 
  resolution of the singularity at the origin of $U_p/G_p$. Repeat 
  this procedure for every singular point $f(p) \in X/G$. We get a 
  resolution of singularities $Y \rightarrow X/G$. Since all the $p$s 
  are fixed points of $T$, $Y \rightarrow X/G$ is $T$-equivariant. The 
  regular map $Y \rightarrow X/G \rightarrow \C^2/G$ gives us our 
  desired $T$-equivariant resolution $\pi: Y\rightarrow \C^2/G$.
  Let $\set{D_i}_{I_Y}$ denote the set of divisors on $Y$ 
  containing the exceptional curves of $\pi$ plus the proper transforms 
  of the divisors $f(D_j)$ for $j\in I_X$. Then $\set{D_i}_{I_Y}$ is a 
  system of $T$-invariant normal crossing divisors on $Y$. 
  
  Since the techniques of section \ref{Poly} pertain to compact varieties, 
  we now describe some natural compactifications of $X$ and $Y$. We 
  may view $\C^2 = \set{(x,y)}$ as sitting inside 
  $\Proj^2=\set{[x:y:z]}$ with the actions of $G$ 
  and $T$ extending to $\Proj^2$ in the obvious manner. Let $\cl{X} 
  \supset X$ be the blow-up of $\Proj^2$ at the origin. Notice that 
  $\cl{X}/G$ is smooth along the hyperplane at infinity. Thus, let 
  $\cl{Y}$ be the resolution of $\cl{X}/G$ which is equal to 
  $\cl{X}/G$ at infinity and which coincides with $Y$ everywhere 
  else. Let $\set{\cl{D}_j}_{I_X}$ and $\set{\cl{D}_i}_{I_Y}$ 
  denote the closures of the corresponding divisors on $X$ and $Y$. Note 
  that the compactified divisors $\set{\cl{D}_j}_{I_X}$ and $\set{\cl{D}_i}_{I_Y}$ 
  do not contain any new points in their intersection loci. Finally 
  let $\Sigma_{\cl{X}}$ and $\Sigma_{\cl{Y}}$ be the polyhedral 
  complexes associated to the above $T$-invariant normal crossing 
  divisors. 
  
  We now describe the relationship between the polyhedral complexes 
  of $\cl{X}$ and $\cl{Y}$. We first set up some notation. As above, let $f:\cl{X}\rightarrow 
  \cl{X}/G$ be the global quotient. Let $p \in \cl{D}_p \cap E$. Then 
  $p$ corresponds to a $2$-dimensional cone $C_p \subset \Sigma_{\cl{X}}$ 
  with lattice $N_p$. Let $C_{g\cdot p}$ be the corresponding $2$-dim 
  cone for each point $g\cdot p$ in the $G$-orbit of $p$. We have an 
  obvious identification $N_{g\cdot p} \cong N_p$. Let 
  $\Sigma_{f(p)}$ be the fan which corresponds to the toric 
  desingularization $\Bl{U_p/G_p}$ of $\C^2/G_p$ described above. Let 
  $N_{f(p)}$ denote the refinement of the lattice $N_p$ which induces 
  the toric quotient map $\C^2\rightarrow \C^2/G_p$. Next, if $D$ is 
  a divisor in $\set{\cl{D}_j}_{I_X}$, let $C_D$ denote the 
  corresponding ray in $\Sigma_{\cl{X}}$, with lattice $N_D \cong \Z$. 
  We define $N_{f(D)}$ to be the suplattice of $N_D$ whose index is 
  equal to the ramification index of $f$ along $D$. Given these data, 
  we obtain the lattice $\Sigma_{\cl{Y}}$ from $\Sigma_{\cl{X}}$ as 
  follows: First identify all the cones $(C_{g\cdot p},N_{p})$ for 
  every point $g\cdot p$ in the $G$-orbit of $p$. Next replace their 
  representative with the toric fan $(\Sigma_{f(p)},N_{f(p)})$. Finally, 
  replace the rays $(C_D,N_D)$ with the rays $(C_{f(D)},N_{f(D)})$. 
  
  Unfortunately, the map $\cl{X}\dashrightarrow \cl{Y}$ is not 
  regular. To remedy this, for every $2$-dimensional cone $C_p \subset 
  \Sigma_{\cl{X}}$, let $\Sigma_p$ be a subdivision of 
  $\Sigma_{f(p)}$ such that the morphism $(\Sigma_p,N_p)\rightarrow 
  (\Sigma_{f(p)},N_{f(p)})$ induces a smooth map of toric varieties. 
  Modify $\Sigma_{\cl{X}}$ by replacing every cone $C_p$ with 
  $\Sigma_p$. Naturally, we identify $\Sigma_{g\cdot p}=\Sigma_{p}$ for 
  every $g\in G$. Define 
  $\hat{X}$ to be the toroidal modification of $\cl{X}$ whose polyhedral 
  complex $\Sigma_{\hat{X}}$ corresponds to the above described 
  modification of $\Sigma_{\cl{X}}$. It is easy to see that 
  $\hat{X}\rightarrow \cl{X}$ factors into a sequence of $G\times T$-equivariant
  blow-ups at complete intersection points. Moreoever, the smooth 
  $T$-map $\mu:\hat{X}\rightarrow \cl{Y}$ is a toroidal morphism. We 
  have the following commutative diagram:
  $$\begin{CD}
     \hat{X}       @>\mu>>     \cl{Y}\\
    @V\phi VV                 @VV\pi V\\
     \Proj^2       @>\psi>>   \Proj^2/G\\ 
\end{CD}$$
 
Let $\set{\hat{D}_j}_{j\in I_{\hat{X}}}$ denote the $T$-invariant 
normal crossing divisors on $\hat{X}$ which correspond to the 
$1$-dimensional rays of $\Sigma_{\hat{X}}$. We have that 
$\phi^*K_{\Proj^2} = K_{\hat{X}}+\sum_{I_{\hat{X}}}\beta_j\hat{D}_j$ for 
some integers $\beta_j < 1$. 

\begin{lem}
    $\phi_*\mathcal{E}ll_{orb}^T(\hat{X},\sum \beta_j\hat{D}_j,G)=
    \mathcal{E}ll_{orb}^T(\Proj^2,0,G)$.
\end{lem}

\begin{proof}
    This follows almost directly from the change of variables formula 
    for the orbifold elliptic genus discussed in the previous 
    section. The only point which requires comment is the fact that 
    the divisors $\set{\phi(\hat{D}_j)}_{I_{\hat{X}}-I_{excep}}$ 
    (where $I_{excep}$ ranges over the exceptional curves) have their 
    base locus at $[0:0:1]$ and therefore do not form a system of 
    normal crossing divisors. However, for any commuting pair $g,h$, 
    we can find a subset of 
    $\set{\phi(\hat{D}_j)}_{I_{\hat{X}}-I_{excep}}$ which consists of 
    normal crossing divisors, and such that every component of 
    $(\Proj^2)^{g,h}$ is inside the intersection locus of this 
    subset. Since the change of variables formula applies to every 
    commuting pair individually, this completes the proof. 
\end{proof}

We also have that $\pi^*K_{\Proj^2/G} = 
K_{\cl{Y}}+\sum_{I_{\cl{Y}}}\alpha_i\cl{D}_i$ for rationals 
$\alpha_i < 1$. Since $\psi^*K_{\Proj^2/G} = K_{\Proj^2}$, we get 
that $\mu^*(K_{\cl{Y}}+\sum_{I_{\cl{Y}}}\alpha_i\cl{D}_i) =
\mu^*\pi^*K_{\Proj^2/G} = \phi^*K_{\Proj^2} = 
K_{\hat{X}}+\sum_{I_{\hat{X}}}\beta_j\hat{D}_j$.

\begin{lem}
    $\mu_*\mathcal{E}ll_{orb}^T(\hat{X},\sum\beta_j\hat{D}_j,G)
    =\big(\frac{2\pi i\theta(-z)}{\theta'(0)}\big )^{2}
    \mathcal{E}ll_T(\cl{Y},\sum\alpha_i\cl{D}_i).$
\end{lem}

\begin{proof}
    Since $\mu : \hat{X}\rightarrow \cl{Y}$ is a toroidal morphism, we 
    have 
    $\mu^*\Omega(\log\sum\alpha_i\cl{D}_i) = 
    \Omega(\log\sum\beta_j\hat{D}_j)$. Therefore, by an argument 
    analogous to the proof of lemma \ref{ChernRoots}:
    $$\frac{c_T(T\hat{X})}{\mu^*c_T(T\cl{Y})}=
    \frac{\prod_{I_{\hat{X}}}(1+c^T_1(\hat{D}_j))}
    {\prod_{I_{\cl{Y}}}(1+\mu^*c_1^T(\cl{D}_i))}.$$
    Following the same argument as in the proof of theorem 
    \ref{OrbChVar}, we 
    are reduced to proving:
    $$\mu_*\frac{1}{|G|}\sum_{gh=hg; {\hat{X}}^{g,h}_\gamma}
    \Phi^T_{\hat{X}^{g,h}_\gamma}\cdot\prod_{I_{\hat{X}}-I^{g,h}_\gamma}
    \orbellipar{\pii{\hat{D}_j}}{\beta_j}{\pii{\hat{D}_j}}\cdot$$
    $$\prod_{I^{g,h}_\gamma}
    \orbnormTh{\pii{\hat{D}_j}+\epsilon_j(g)-\epsilon_j(h)\tau}{\beta_j}
    e^{2\pi i(-\beta_j+1)\epsilon_j(h)z}=$$
    
    $$\prod_{I_{\cl{Y}}}\orbellipar{\pii{\cl{D}_i}}{\alpha_i}{\pii{\cl{D}_i}}$$
    
    Let $H \in \C[[\Sigma_{\cl{Y}}]]$ be the piecewise convergent 
    power series:
    
    $$H|_C = \prod_{i=1}^{\dim C}\orbellipar{\pii{x^C_i}}{\alpha_i}
    {\pii{x^C_i}}.$$
    
    Similarly, let $F^{g,h}_\gamma \in \C[[\Sigma_{\hat{X}}]]$ be 
    defined as in the proof of theorem \ref{OrbChVar}. By lemma 
    \ref{PushForward}, we are 
    reduced to proving:
    $$\nu_*\frac{1}{|G|}\sum_{gh=hg; \gamma}F^{g,h}_\gamma = H.$$
    
    Let $C \subset \Sigma_{\cl{Y}}$ be a $2$-dimensional cone with 
    lattice $N_C$. Let 
    $\nu^{-1}C$ be the collection of fans $\set{\Sigma_C^k}_k$. From 
    our construction of $\Sigma_{\hat{X}}$ and $\Sigma_{\cl{Y}}$, the 
    $\Sigma_C^k$s are isomorphic to a fixed subdivision $\Sigma_C$ 
    with lattice $N(\Sigma_C)$. Let $G = N_C/N(\Sigma_C)$ and $d$ be 
    the cardinality of the set $\set{\Sigma_C^k}$. We clearly have 
    $d\cdot|G(\Sigma_C)| = |G|$. By the pushforward formula for 
    toroidal morphisms:
    $$\bigg\{\nu_*\frac{1}{|G|}\sum_{gh=hg;\gamma}F^{g,h}_\gamma
    \bigg\}_C=
    \frac{d}{|G|}\sum_{C_i\subset \Sigma_C}\sum_{gh=hg;\gamma}
    (F^{g,h}_\gamma)_{C_i}\frac{\prod_{j=1}^2 x^C_j}
    {\prod_{j=1}^2 x^{C_i}_j}=$$ 
    
    $$\frac{1}{|G(\Sigma_C)|}\sum_{C_i\subset \Sigma_C}
    \sum_{g,h\in G(\Sigma_C)}$$
    $$\prod_{j=1}^2
    \orbnormTh{\pii{x^{C_i}_j}+\epsilon_j(g)-\epsilon_j(h)\tau}
    {\beta_j}e^{2\pi i(-\beta_j+1)\epsilon_j(h)z}\cdot\prod_{j=1}^2x^C_j.$$
    We wish to show that this last expression is equal to
    $$\prod_{j=1}^2\orbellipar{\pii{x^C_j}}{\alpha_j}{\pii{x^C_j}}.$$
    Since the coefficients $\set{-\beta_j+1}$ are the pullbacks of the 
    coefficients $\set{-\alpha_i+1}$, this identity follows from 
    Lemma $8.1$ in \cite{BL}. This completes the proof of the lemma.
\end{proof}

Now let $p: V\rightarrow \C^2/G$ be the $T$-equivariant crepant 
resolution, i.e., $K_V = p^*K_{\C^2/G}$. By the equivariant 
factorization theorem for surfaces, we may connect $Y$ to $V$ by a 
finite sequence of equivariant blow-ups and blow-downs. In other 
words, we may form the commutative diagram:
$$\begin{CD}
     \Bl{Y}        @>f>>            Y\\
    @VgVV                     @VV\pi V\\
     V             @>p>>         \C^2/G\\ 
\end{CD}$$

Here $f$ and $g$ are sequences of $T$-equivariant blow-ups. Moreover, 
we may assume that the blow-ups of $f$ occur at complete intersection 
points. It is not immediately clear that the same holds true for $g$. 
However, every blow-up point for $g$ will occur at a $T$-fixed point 
$p \in D$ for some $D \cong \Proj^1$, with $D$ $T$-invariant. From 
toric geometry, every such $p$ may be represented as a complete 
intersection of $T$-invariant divisors in a neighborhood of $D$. By 
functorial localization, we may therefore assume that $p$ is a 
complete intersection point when we calculate pushforwards.

Define $\Bl{D}_i$ and $\tilde{\alpha}_i$ so that 
$f^*(K_Y+\sum\alpha_iD_i)= K_{\Bl{Y}}+\sum\tilde{\alpha}_i\Bl{D}_i$. 
Then $K_{\Bl{Y}}+\sum\tilde{\alpha}_i\Bl{D}_i = 
f^*\pi^*K_{\C^2/G}=g^*K_V$. By the equivariant change of variables 
formula and functorial localization:
$$\sum_{P\subset V}\int_P\frac{\mathcal{E}ll_T(V,0)}{e(P)} =
\sum_{F\subset Y}\int_F\frac{\mathcal{E}ll_T(Y,\sum\alpha_iD_i)}
{e(F)}.$$
Here $\set{P}$ and $\set{F}$ run over the fixed components of $V$ 
and $Y$, respectively. The above formula combined with the previous 
two lemmas imply the following equivariant elliptic genus analogue of 
the classical McKay correspondence:

\begin{thm}
    $Ell^T_{orb}(\C^2,0,G)=\big(\frac{2\pi i\theta(-z)}{\theta'(0)}\big )^{2}
    Ell_T(V,0).$
\end{thm}

Letting the parameter $z\to 0$ we get the following corollary:

\begin{cor}
    The number of irreducible representations of $G$ is equal to 
    $\sum_{P\subset V}\int_P e(TP)$.
\end{cor}

In other words, $e_{orb}(\C^2,G) = e(V)$.

\section{Relation to Batyrev's Stringy Euler Number}\label{StringyEuler}
Let $X$ be a smooth quasiprojective variety and $D = \sum_I a_iD_i$ 
an effective divisor with simple normal crossings. In \cite{B} Batyrev 
defines the stringy Euler number of the pair $(X,D)$ as:
$$e_{str}(X,D) = \sum_{J\subset I}e(D_J^o)\prod_{j\in J}
\frac{1}{a_j+1}.$$
Here $D_J^o = \cap_{j\in J}D_j - \cup_{I-J}D_i$. The definition is 
best understood from the point of view of motivic integration. 

Suppose that $X$ has a $T$ action with compact fixed components, and 
that the irreducible components of $D$ are $T$-invariant. Then the 
equivariant elliptic genus $Ell_T(X,-D)$ is well-defined and it is 
natural to question how it is related to Batyrev's stringy Euler 
number. In this section we prove the following relation:

\begin{thm}
Let $X$ and $D$ be as above. Then:
$$\lim_{y\to 1}\lim_{q\to 0}
\bigg(\frac{2\pi i\theta(-z)}{\theta'(0)}\bigg )^{\dim X}
Ell_T(X,-D) = e_{str}(X,D).$$
\end{thm}

\begin{proof}
We first note that we may rewrite $e_{str}(X,D)$ in the more 
convenient form:
$$e_{str}(X,D) = \sum_{J\subset I}e(D_J)\prod_{j\in 
J}(\frac{1}{a_j+1}-1).$$
Here $D_J = \cap_{j\in J}D_j$. For a quasiprojective $T$-space $V$ 
with compact fixed components, let $\chi_{-y}^T(V)$ denote the 
equivariant $\chi_{-y}$-genus of $V$; that is, 
$$\chi_{-y}^T(V) = \sum_{F\subset V^T}\int_F \prod_{TF}
\frac{f_i(1-ye^{-f_i})}{1-e^{-f_i}}\prod_{\nu_F}\frac{1-ye^{-n_j-w_j}}{1-e^{-n_j-w_j}}.$$
Here $n_j$ are the formal chern roots of the normal bundle $\nu_F$ to 
$F$ and $w_j$ are the infinitesimal weights of the $T$ action on the 
fibers of $\nu_F$. It is easy to verify that:
$$\lim_{y\to 1}\sum_{J\subset I}\chi_{-y}^T(D_J)\prod_{j\in J}
(\frac{y-1}{y^{a_j+1}-1}-1) = e_{str}(X,D).$$
The above equality essentially follows from the fact that $e(V)$ 
corresponds to  
$\sum_{F \subset V^T}e(F)$, where $\set{F}$ denotes the collection of 
fixed components of $V$ (all of which are compact). Thus, we are 
reduced to identifying $\sum_{J\subset I}\chi_{-y}^T(D_J)\prod_{j\in 
J}(\frac{y-1}{y^{a_j+1}-1}-1)$ with a specialization of the 
equivariant elliptic genus.

Fix an indexing set $J \subset I$ and a fixed component $F$ of $X$ 
which intersects $D_J$ nontrivially. Write $\set{D_j}_{j\in J}$ as 
the disjoint union of collections $\set{D_\alpha} \cup 
\set{D_\beta}$, where $D_\beta$ denote the divisors which contain 
$F$. Note that since $F \cap D_J \neq \emptyset$, for any of the 
divisors $D_\alpha$, $D_\alpha \cap F$ is a proper subset of $F$. 

The collection $\set{D_J \cap F_i}_{F_i \in \hbox{Fix}(X)}$ describes 
a partition of the collection of fixed components of $D_J$. Let $P = 
D_J \cap F$. From the above discussion, $\nu_{P/F} = 
\oplus_{\alpha}\mathcal{O}(D_\alpha)|_P$.

Consider the integral:
$$\int_F\prod_{TF}\frac{f_i(1-ye^{-f_i})}{1-e^{-f_i}}
\prod_{\nu_{F/X}}\frac{1-ye^{-n_\ell-w_\ell}}{1-e^{-n_\ell-w_\ell}}
\prod_{\beta}\frac{1-e^{-D_\beta-w_\beta}}{1-ye^{-D_\beta-w_\beta}}
\prod_\alpha \frac{1-e^{-D_\alpha}}{1-ye^{-D_\alpha}}.$$
Letting $N = \nu_{F/X}-\oplus_\beta\mathcal{O}(D_\beta)$,
we may rewrite the integral as:
$$\int_F\prod_{TF}\frac{f_i(1-ye^{-f_i})}{1-e^{-f_i}}
\prod_{N}\frac{1-ye^{-n_\ell-w_\ell}}{1-e^{-n_\ell-w_\ell}}
\prod_\alpha 
\frac{1-e^{-D_\alpha}}{D_\alpha(1-ye^{-D_\alpha})}\prod_\alpha 
D_\alpha = $$
$$\int_P\prod_{TP}\frac{p_k(1-ye^{-p_k})}{1-e^{-p_k}}
\prod_{\nu_{P/D_J}}\frac{1-ye^{-n_\ell-w_\ell}}{1-e^{-n_\ell-w_\ell}}.$$
Summing over all $J \subset I$, we therefore have:
\begin{align*}
{}&\sum_{J\subset I}\chi_{-y}^T(D_J)\prod_{j\in J}(\frac{y-1}{y^{a_j+1}-1}-1) = \\
{}&\sum_{J,F} \int_F\prod_{TF}\frac{f_i(1-ye^{-f_i})}{1-e^{-f_i}}
\prod_{\nu_{F/X}}\frac{1-ye^{-n_\ell-w_\ell}}{1-e^{-n_\ell-w_\ell}}
\prod_J\frac{1-e^{-D_j-w_j}}{1-ye^{-D_j-w_j}}(\frac{y-y^{a_j+1}}{y^{a_j+1}-1}) \\ 
{}&=\sum_F\int_F\prod_{TF}\frac{f_i(1-ye^{-f_i})}{1-e^{-f_i}}
\prod_{\nu_{F/X}}\frac{1-ye^{-n_\ell-w_\ell}}{1-e^{-n_\ell-w_\ell}}
\prod_I\frac{1-e^{-D_i-w_i}}{1-ye^{-D_i-w_i}}
\frac{y-1}{y^{a_i+1}-1} \\
\end{align*}
It is easy to see that this last expression is equal to \newline
$\lim_{q\to 0}\bigg(\frac{2\pi i\theta(-z)}{\theta'(0)}\bigg )^{\dim X}Ell_T(X,-D)$. 
This completes the proof.
\end{proof}

\section{Appendix}
\begin{lem}
   Let $f:X \rightarrow Y$ be a $T$-map of smooth compact simply 
    connected complex manifolds. 
    Let $D \subset Y$ be a $T$-invariant divisor and let $E_i$ be 
    $T$-invariant normal crossing divisors on $X$ such that $f^*D = 
    \sum a_i E_i$ as Cartier divisors. Then for any 
    $\varepsilon$-regular neighborhood $U_{\varepsilon}$ of $D$ 
    there exist generators $\Theta_{E_i}^T$ for $c_1^T(E_i)$ and 
    $\Theta_D^T$ for $c_1^T(D)$ with the following properties:
    
    $(1)$ $\Theta_D^T$ has compact support in $U_{\varepsilon}$ and 
    $\Theta_{E_i}^T$ have compact support in $f^{-1}(U_{\varepsilon})$.
    
    $(2)$ $f^* \Theta_D^T = \sum a_i \Theta_{E_i}^T + d_T (\eta)$ on 
    the level of forms, where $\eta$ is a $T$-invariant form with 
    compact support in $U_{\varepsilon}$.
    
    $(3)$ $\Theta_D^T$ and $\Theta_{E_i}^T$ represent the extension by zero of the equivariant 
    Thom classes of the varieties $D$ and $E_i$
\end{lem}

The only real issue above is to ensure that $\eta$ has compact support 
in the desired neighborhood.

\begin{proof}
    We first solve this problem in the non-equivariant category. For 
    $V$ any Cartier divisor, denote by $L_V$ the line bundle it 
    induces. Let $U_{\varepsilon}$
    be a $T$-invariant tubular neighborhood of $D$ of radius 
    $\varepsilon$. Outside $U_{\frac{\varepsilon}{2}}$, the constant 
    function $1$ is a section of $L_D$. Define a metric $h_{far}$ in this 
    region by $h_{far} = \norm{1}^2 \equiv 1$. Let $h_{near}$ be a metric 
    inside $U_{\varepsilon}$. Piece the two metrics into a global 
    metric $h$ on $L_D$ using a partition of unity. The first chern class 
    of $L_D$ is then represented by the form 
    $\Theta_D =\frac{i}{2\pi}\delbar\del \log h$. This form clearly has compact 
    support in $U_{\varepsilon}$.
    
    Let $U_{\varepsilon_i}$ be tubular neighborhoods of $E_i$. 
    Choose $\varepsilon_i$ small enough so that each of these 
    neighborhoods is contained in $f^{-1}U_{\varepsilon}$. Define 
    metrics $h_i$ on $E_i$ in a manner analogous to the above 
    construction of $h$. Clearly the forms $\Theta_{E_i} = 
    \frac{i}{2\pi}\delbar\del \log h_i$ have compact support in 
    $U_{\varepsilon_i}$ and represent the first chern classes of $L_{E_i}$. 
    
    We have two natural choices for a metric on $f^*L_D$, namely 
    $f^*h$ and $h_1^{a_1}\cdots h_k^{a_k}$. Choose a smooth nonzero 
    function $\varphi$ so that $f^*h = \varphi h_1^{a_1}\cdots h_k^{a_k}$.
    Notice that $\varphi \equiv 1$ outside $f^{-1}U_{\varepsilon}$. 
    We have:
    $$\delbar\del\log f^*h = \delbar\del\log\varphi +
    \sum_i a_i\delbar\del\log h_i.$$
    But this implies that $f^*\Theta_D = \sum_i a_i\Theta_{E_i} 
    +\frac{i}{2\pi}\delbar\del\log\varphi$. If we let $d^c = 
    \frac{i}{4\pi}(\delbar-\del)$, we may write this last equation as:
    $$f^*\Theta_D = \sum_i a_i\Theta_{E_i}-dd^c\log\varphi.$$
    
    The form $\eta = -d^c\log\varphi$ clearly has compact support in 
    $f^{-1}U_{\varepsilon}$. It remains to argue that $\Theta_D$ and 
    $\Theta_{E_i}$ represent the Thom classes of $D$ and $E_i$. It is a standard fact that these classes are Poincar\'e duals to their respective divisors. If a divisor is homologously non-trivial, then clearly its Thom classes coincides with its Poincar\'e dual. If a divisor is homologously trivial, then it must follow that the extension by zero of its Thom class is trivial. Either way this implies that the above classes represent the extension by zero of the Thom classes of their respective divisors. This completes the non-equivariant portion of the proof.    
    
    By averaging over the group $T$, we may assume that all the forms 
    above are $T$-invariant. For notational simplicity, let us assume 
    that $T = S^1$. Let $V$ be the vectorfield on $X$ induced by the 
    $T$-action. Let $g_i$ be the functions compactly supported in 
    $f^{-1}U_{\varepsilon}$ which satisfy the moment map equation 
    $i_V\Theta_{E_i} = dg_i$. Similarly, let $W$ be the vectorfield 
    on $Y$ defined by the $T$-action and define $g$ so that it 
    satisfies $i_V\Theta_D = dg$ and has support inside 
    $U_{\varepsilon}$. Note that since $f$ is $T$-equivariant, 
    $i_Vf^*\Theta_D = f^*i_W\Theta_D = f^*dg$. We 
    then have $d(g\circ f) = \sum_i a_i dg_i + i_Vd\eta = \sum_i a_i dg_i 
    -di_V\eta$. Hence $g\circ f = \sum_i a_i g_i -i_V\eta$. But this implies 
    that:
    $$f^*(\Theta_D+g) = \sum_i a_i(\Theta_{E_i}+g_i)+(d-i_V)\eta.$$
    But this is precisely the relation we wish to in the equivariant 
    cohomology.
\end{proof}

\begin{lem} \bf{(Excess Intersection Formula)} \it
    Let $X$ be a smooth compact variety with irreducible normal crossing 
    divisors $D_1,\ldots, D_k$. For $I \subset \set{1,\ldots,k}$ 
    denote by $X_{I,j}$ the $j$th connected component of 
    $\cap_{I}D_i$ and by $\Phi_{I,j}$ its Thom class. Fix irreducible 
    subvarieties $X_{I_1,j_1}$ and $X_{I_2,j_2}$. For $I_0 = I_1 \cup 
    I_2$, let $X_{I_0,j}$ be the irreducible components of 
    $X_{I_1,j_1} \cap X_{I_2,j_2}$. Then:
    $$\Phi_{I_1,j_1}\wedge \Phi_{I_2,j_2} = 
    \sum_{I_0,j}\Phi_{I_0,j}\prod_{I_1\cap I_2}\Phi_{i}.$$
\end{lem}

\begin{proof}
    Let $N_{I,j}$ be tubular neighborhoods of $X_{I,j}$ which are 
    disjoint for each indexing set $I$ and which satisfy 
    $N_{I,j} \subset N_{I',j'}$ for $X_{I,j}\subset X_{I',j'}$. If we choose 
    $\Phi_i$ to have compact support in a sufficiently small tubular 
    neighborhood of $D_i$, then $\prod_I \Phi_i$ will have compact 
    support in $\coprod_{j}N_{I,j}$. Moreover, the extension by zero 
    of $(\prod_I \Phi_i)|_{N_{I,j}}$ will represent the Thom class of 
    $X_{I,j}$ (see \cite{BT}). We may also ensure that 
    $\Phi_{I_1,j_1}\wedge \Phi_{I_2,j_2}$ has compact support in 
    $\coprod_j N_{I_0,j}$. Thus:
    $$\Phi_{I_1,j_1}\wedge \Phi_{I_2,j_2} = \sum_{I_0,j}
    \big( \prod_{I_1}\Phi_i \prod_{I_2}\Phi_i \big)|_{N_{I_0,j}} =
    \sum_{I_0,j}\big(\prod_{I_0}\Phi_i \prod_{I_1\cap 
    I_2}\Phi_i\big)|_{N_{I_0,j}} = $$
    $$\sum_{I_0,j}\big(\prod_{I_0}\Phi_i \big)|_{N_{I_0,j}} \prod_{I_1\cap 
    I_2}\Phi_i.$$
    This yields the desired formula.
\end{proof}

\begin{rmk}
    \rm Note that the above proof clearly extends to the equivariant 
    category.
\end{rmk}

\end{document}